\documentclass[3p]{elsarticle}

\usepackage{amssymb}
\usepackage{amsmath}
\usepackage{tabu}
\usepackage[tableposition=top]{caption}
\usepackage{booktabs}
\usepackage[utf8]{inputenc}
\usepackage{epsfig,graphicx}
\usepackage{epstopdf}
\usepackage{graphicx}
\usepackage{caption}
\usepackage{subcaption}
\usepackage{hyperref}

\usepackage[dvipsnames,svgnames,x11names]{xcolor}

\usepackage[final]{changes}

\definechangesauthor[color=red]{Rev.1}
\definechangesauthor[color=blue]{Rev.3}
\definechangesauthor[color=green]{Authors}
\definechangesauthor[color=teal]{Rev.2}
\definechangesauthor[color=orange]{Revs}
\usepackage{url}
\hypersetup{
    colorlinks,
    citecolor=black,
    filecolor=black,
    linkcolor=black,
    urlcolor=black
}

\biboptions{numbers,sort&compress}

\journal{Journal of Computational Physics}

\makeatletter
\def\ps@pprintTitle{%
 \let\@oddhead\@empty
 \let\@evenhead\@empty
 \def\@oddfoot{}%
 \let\@evenfoot\@oddfoot}
\makeatother

\begin{document}

\begin{frontmatter}

\title{An acoustic-convective splitting-based approach for the Kapila two-phase flow model}

 \author[label1,label2]{M.F.P. ten Eikelder\corref{cor1}}
 \cortext[cor1]{Corresponding author. \\ \noindent Present address: Delft University of Technology, Department of Mechanical, Maritime and Materials Engineering, P.O. Box 5, 2600 AA Delft, The Netherlands}
 \ead{m.f.p.teneikelder@tudelft.nl}
 \author[label1,label3]{F. Daude}
 \author[label2]{B. Koren}
 \author[label2]{A.S. Tijsseling}
 \address[label1]{EDF R\&D, AMA, 7 boulevard Gaspard Monge 91120, Palaiseau, France}
 \address[label2]{Eindhoven University of Technology, Department of Mathematics and Computer Science, P.O. Box 513, 5600 MB Eindhoven, The Netherlands}
 \address[label3]{IMSIA, UMR EDF-CNRS-CEA-ENSTA 9219, Universit\'{e} Paris Saclay, 828 Boulevard des Mar\'{e}chaux 91762 Palaiseau, France}

\begin{abstract}
In this paper we propose a new acoustic-convective splitting-based numerical scheme for the Kapila five-equation two-phase flow model. The splitting operator decouples the acoustic waves and convective waves. The resulting two submodels are alternately numerically solved to approximate the solution of the entire model. The Lagrangian form of the acoustic submodel is numerically solved using an HLLC-type Riemann solver whereas the convective part is approximated with an upwind scheme. The result is a simple method which allows for a general equation of state. Numerical computations are performed for standard two-phase shock tube problems. A comparison is made with a non-splitting approach. The results are in good agreement with reference results and exact solutions.
\end{abstract}


\end{frontmatter}

%
%

\section{Introduction}\label{section Introduction}
Compressible two-phase and two-fluid flow phenomena arise in many natural features and industrial applications. Examples are groundwater flow, surface wave impacts, oil slicks, water-air flows, shock-bubble interaction and (condensation induced) water hammer phenomena. The study of two-phase flow is a challenging research area which is of interest to both engineers and scientists.\\
\indent Various models can be used to describe two-phase flows. Many of these models can be classified as two-fluid models, or homogeneous models. Among the two-fluid flow models, which are generally considered as the most complete, the model of Baer and Nunziato \cite{baer1986two} is one of the best known. This model consists of equations for each of the two fluids' mass, momentum, energy, and of an equation describing the topology of the two-fluid interface. Romenski et al. \cite{romenski2007conservative} proposed a seven-equation model for two-phase compressible flow which can be written in Baer-Nunziato form in the heat flux relaxation limit. Due to the complexity of the seven-equation models, linked to their large number of different waves \cite{crouzet2013approximate, saurel1999multiphase, massoni2002proposition, schwendeman2006riemann, gallouet2004numerical, tokareva2010hllc, dumbser2011simple, ambroso2012godunov, herard2012fractional, crouzet2015validation, lochon2016comparison, daude2016computation}, reduced models with less equations have been proposed.\\
\indent The five-equation models form an important class of reduced models. The original five-equation two-phase flow model of Kapila et al. \cite{kapila2001two} has been derived from the two-fluid flow model of Baer and Nunziato. To study pure interface problems the model of Allaire et al. \cite{allaire2002five, daude2014numerical, kokh2010anti} can be used. The model of Kapila et al., describing inviscid, non-heat-conducting, compressible two-fluid flow, allows for mixtures. To model phase transitions, the five-equation model has been extended by taking temperature and chemical potential relaxation effects into account \cite{saurel2008modelling}. Murrone and Guillard \cite{murrone2005five} give an analysis of the five-equation model and indicate that the five-equation model is a good approximation of the seven-equation two-fluid model. Kreeft and Koren \cite{kreeft2010new} propose a new formulation of the five-equation model, in which the topological equation is replaced by an energy equation. An Osher-type approximation is used for the evaluation of the fluxes and the energy\added[id=Authors]{-}exchange term in the discretized system. Ahmed et al. \cite{ahmed2015central} use a central upwind scheme for the new formulation to study shock-bubble interaction problems. Daude et al. \cite{daude2014numerical} present computations with the original five-equation model of Kapila et al. using an HLLC-type scheme in the context of an Arbitrary Lagrangian-Eulerian formulation.\\
\indent Serious difficulties are posed by the non-conservative terms in the topology equation of the five-equation model. In particular, (i) approximating the term containing velocity divergence, (ii) performing shock computations with a non-conservative model \added[id=Authors]{and (iii) ensuring volume fraction positivity} \cite{berry2009simple, petitpas2007relaxation, saurel2007relaxation} is difficult. Several approaches have been suggested to circumvent these issues. Abgrall and Perrier \cite{abgrall2006asymptotic} present, using probabilistic multiscale interpretation of multiphase flows, a locally conservative scheme to tackle the issues. Saurel-Petitpas-Berry \cite{berry2009simple, saurel2009simple} propose to relax the pressure equilibrium assumption and obtain a non-conservative hyperbolic six-equation model which simplifies numerical resolution. Jiang et al. \cite{Jiang} use this six-equation approach with a novel mass transfer between liquid and vapor.\\
\indent The aim of the present paper is to propose an acoustic-convective splitting-based numerical method for the five-equation two-phase flow model. Due to its simplicity, the original five-equation model of Kapila et al., without any relaxation or modification, is considered. \added[id=Authors]{Furthermore, the speed of sound of this model corresponds to the Wood speed of sound which is known to be in good agreement with the experimental data obtained at moderate frequencies of sound (pressure disturbance) in air-water mixtures.} The present approach is inspired by the Lagrange-Projection-like scheme originally proposed for the Euler equations of gas dynamics, by Chalons et al. \cite{chalons2014all}. \added[id=Authors]{In this paper a method similar to that} from \cite{chalons2014all} is extended to the full two-phase five-equation model. Related work of the authors about the splitting approach has been presented in \cite{ten2016lagrange}. \added[id=Authors]{Our scheme uses an HLLC-type scheme for the acoustic model and a classical upwind scheme for the convective model.} \added[id=Rev.3]{Conservation of mass, momentum, energy and partial mass, as well as the positivity of the volume fraction and the mass fraction are ensured.} \added[id=Revs]{The advantages of the proposed approach are (i) its simplicity and (ii) its accurate capturing at shock waves and (iii) the potential to deal with low-Mach number flows. Approximate Godunov approaches and direct approaches may lead to inaccuracies at highly subsonic flows. By using a splitting operator these inaccuracies can be prevented \cite{chalons2014all}. Furthermore, unlike Osher-type schemes \cite{kreeft2010new}, the current approach can deal with a general equation of state (just like the direct approach from \cite{daude2014numerical}).} A similar idea \added[id=Authors]{has been} proposed by Huber et al. \cite{huber2015time}. They use a compressible projection method with a level-set method describing the interface motion to study the interaction of an ultrasound wave with a bubble. \\
\indent The paper is organized as follows. In Section \ref{section Two-phase flow model} the five-equation two-phase Kapila et al. flow model is shortly rehearsed. The novel acoustic-convective splitting scheme is presented in Section \ref{section Numerical scheme}. The numerical scheme is assessed for shock-tube problems in Section \ref{section Numerical results}\added[id=Revs]{, and a comparison with the direct approach is made in terms of accuracy, efficiency and robustness.} Conclusions are drawn in Section \ref{section Conclusions}.

%
%

\section{Two-phase flow model}\label{section Two-phase flow model}
The five-equation model of Kapila et al. \cite{kapila2001two} describes the dynamics of inviscid two-phase flows evolving in mechanical equilibrium (i.e. equilibrium of velocity and pressure is assumed across the fluid interface). The model consists of four balance equations for conservative quantities: two for mass (bulk mass and mass of one of the two phases), one for the bulk momentum and one for the bulk total energy. The fifth equation is a topological equation, of non-conservative type, which describes the evolution of the volume fraction. In one dimension, the governing equations read:
  \begin{subequations}\label{governing equations 5 equation model}
    \begin{alignat}{6}
            &\partial_t\rho&+&\partial_x\left( \rho u\right)&&&=0,\label{model bulk mass eq}\\
            &\partial_t(\rho u)&+&\partial_x \left( \rho u^2+p\right)&&&=0,\label{model bulk momentum eq}\\
            &\partial_t(\rho E)&+&\partial_x \left(\rho E u+pu\right)&&&=0,\label{model bulk energy eq}\\
            &\partial_t(\alpha_1 \rho_1)&+&\partial_x \left(\alpha_1\rho_1 u\right)&&&=0,\label{model single mass eq}\\
            &\partial_t\alpha_1&+&u \partial_x \alpha_1&+&K \partial_x u &=0, \label{model volume fraction eq}
        \end{alignat}
    \end{subequations}
where $t$ is the time, $x$ the spatial coordinate, $\rho$ the mixture density, $u$  the bulk velocity, $p$ the pressure and $E$ the mixture total specific energy. The interfacial variable $K$ is specified later. The variable $\alpha_k, k=1,2$, represents the volume fraction of phase $k$, with the saturation constraint $\alpha_1+\alpha_2=1$, and $\rho_k$ denotes the density of \added[id=Authors]{phase} $k$. In terms of separated fluid variables, the bulk density is given by
    \begin{equation}
            \rho=\alpha_1\rho_1+\alpha_2\rho_2.
    \end{equation}
We define the mass fraction $Y_k$ of phase $k$ as $\rho Y_k=\alpha_k \rho_k$. The entropy equations, i.e.:
    \begin{equation}\label{entropy equations}
        \partial_t(\alpha_k \rho_k s_k)+\partial_x(\alpha_k \rho_k s_k u)=0,
    \end{equation}
with $s_k$ the specific entropy of \added[id=Authors]{phase} $k$, complement the model in absence of shocks \cite{murrone2005five}. All the dissipative effects are neglected (inviscid, non-heat conducting flow is considered) and thus it can be written as
    \begin{equation}
            \dfrac{{\rm D} s_k}{{\rm D} t}=0,
    \end{equation}
with the Lagrangian derivative ${\rm D} /{\rm D} t:=\partial_t   + u \partial_x $. The total specific energy of the mixture is given by:
    \begin{equation}
            \rho E=\alpha_1\rho_1E_1+\alpha_2\rho_2E_2,
    \end{equation}
where the total specific energy of each of the two \added[id=Authors]{phases} is
    \begin{equation}
        E_k=e_k+\frac{1}{2}u^2,
    \end{equation}
with $e_k$ the internal specific energy of \added[id=Authors]{phase} $k$. The bulk internal specific energy is given by
    \begin{equation}
        \rho e=\alpha_1\rho_1e_1+\alpha_2\rho_2e_2,
    \end{equation}
and hence,
    \begin{equation}
        E=e+\frac{1}{2}u^2.
    \end{equation}
In the present paper, the model is completed with the stiffened gas (SG) equation of state (EOS) for each phase:
    \begin{equation}\label{SG EOS}
        p=\rho_k\added[id=Authors]{(e_k-\eta_k)}(\gamma_k-1)-\gamma_k\pi_k,
    \end{equation}
where the pressure equilibrium across the interface is used. The ratio of specific heats $\gamma_k$, stiffness $\pi_k$ \added[id=Authors]{and energies at a reference state $\eta_k$} are characteristic constants of the thermodynamic behavior of fluid $k$. Expression (\ref{SG EOS}) reduces to the perfect gas (PG) EOS when $\pi_k$ and $\eta_k$ is equal to zero whereas a large value of $\pi_k$ implies a near-incompressible behavior \cite{flaatten2011solutions}. The SG EOS parameters are determined by shock wave Hugoniot curves \cite{saurel1999simple, coralic2013shock, gojani2009shock}. This EOS is often used as a reasonable approximation for both liquids and gases under high pressure conditions \cite{crouzet2015validation, lochon2016comparison, daude2014numerical, kreeft2010new, abgrall2003discrete, lund2013splitting}.
The EOS allows the determination of the speed of sound of each single phase
    \begin{equation}\label{SG EOS each single phase}
        c_k^2 \equiv \dfrac{p-\rho_k^2\partial_{\rho_k}e_k}{\rho_k^2\partial_pe_k}=\gamma_k\dfrac{p+\pi_k}{\rho_k}.
    \end{equation}
The interfacial variable in the topology equation (\ref{model volume fraction eq}) is given by
    \begin{equation}\label{K}
        K= \left(\rho_1 c_1^2-\rho_2 c_2^2\right)/\left(\frac{\rho_1 c_1^2}{\alpha_1}+\frac{\rho_2 c_2^2}{\alpha_2}\right).
    \end{equation}
The internal specific energy of the mixture satisfies
    \begin{equation}
        \rho e=p\left(\dfrac{\alpha_1}{\gamma_1-1}+\dfrac{\alpha_2}{\gamma_2-1}\right)+\alpha_1\left(\dfrac{\gamma_1}{\gamma_1-1}\pi_1\added[id=Authors]{+\rho_1 \eta_1}\right)+\alpha_2\left(\dfrac{\gamma_2}{\gamma_2-1}\pi_2\added[id=Authors]{+\rho_2\eta_2}\right).
    \end{equation}
The five-equation model (\ref{governing equations 5 equation model}) is hyperbolic and admits the wave speeds \cite{murrone2005five}
\begin{equation}\label{wave speeds}
    \lambda_1=u-c,~~\lambda_{2,3,4}=u,~~\lambda_5=u+c,
\end{equation}
with $c$ the mixture speed of sound which obeys the Wood formula \cite{wood1930textbook}:
    \begin{equation}\label{speed of sound a la Wood}
        \frac{1}{\rho c^2}=\frac{\alpha_1}{\rho_1 c_1^2}+\frac{\alpha_2}{\rho_2 c_2^2}.
    \end{equation}
The characteristic fields associated with the eigenvalues $\lambda_{2,3,4}$ are linearly degenerate (LD) and the other two fields are genuinely nonlinear (GNL) \cite{murrone2005five}.

%
%

\section{Numerical scheme}\label{section Numerical scheme}
A novel splitting-based numerical scheme is presented, leading to two operators: one associated with the pressure and the other \added[id=Authors]{with} the advection. The two submodels are referred to \added[id=Authors]{as} acoustic and convective, respectively, in the sequel. First, the treatment of the acoustic submodel is discussed for which a simple and robust HLLC-type Riemann solver is used. Next, the upwind scheme for the convective submodel is given.


\subsection{The splitting approach}\label{subsection Splitting}
The five-equation model deals with two kinds of wave speeds associated with its eigenvalues, i.e. the GNL waves are linked to acoustic pressure waves whereas the LD wave is connected to the convective velocity. In certain situations such as subsonic flows, the ratio between these two speeds can be large, leading to inaccuracies when using approximate Godunov approaches. In order to decouple acoustic and convective phenomena, a splitting operator is proposed. This splitting is inspired by the one proposed by Chalons et al. \cite{chalons2014all} for the Euler equations of gas dynamics.\\
\indent By using product-rule arguments the Kapila five-equation model (\ref{governing equations 5 equation model}) is split into (i) the acoustic system:
\begin{subequations}\label{acoustic}
    \begin{alignat}{6}
        &\partial_t\rho&+&\rho \partial_xu&&&=0,\label{acoustic rho}\\
        &\partial_t(\rho u)&+&\rho u \partial_xu&+&\partial_xp&=0,\label{acoustic rho u}\\
        &\partial_t(\rho E)&+&\rho E \partial_xu&+&\partial_x(p u)&=0,\label{acoustic rho E}\\
        &\partial_t Y_1 &&&&&=0,\label{acoustic alpha1 rho1}\\
        &\partial_t\alpha_1&+&K \partial_xu&&&=0,\label{acoustic alpha1 rho1 E1}
    \end{alignat}
\end{subequations}
and (ii) the convective system:
    \begin{subequations}\label{convective}
        \begin{alignat}{5}
            &\partial_t\rho&+&u\partial_x \rho&=0,\label{convective rho}\\
            &\partial_t(\rho u)&+&u\partial_x (\rho u)&=0,\label{convective rho u}\\
            &\partial_t(\rho E)&+&u\partial_x(\rho E)&=0,\label{convective rho E}\\
            &\partial_tY_1 &+&u\partial_xY_1&=0,\label{convective alpha1 rho1}\\
            &\partial_t\alpha_1 &+&u\partial_x\alpha_1&=0,\label{convective alpha1 rho1 E1}
        \end{alignat}
    \end{subequations}
where the evolution of the mass fraction,  Eqs. (\ref{acoustic alpha1 rho1}) and (\ref{convective alpha1 rho1}), follows from Eqs. (\ref{model bulk mass eq}) and (\ref{model single mass eq}). The corresponding entropy equations of the acoustic and convective systems are respectively:
\begin{subequations}\label{entropy subsystems}
    \begin{alignat}{6}
        & \partial_t s_k&&&=0,&\label{entropy subsystems: acoustic}\\
        &\partial_t s_k &+&u\partial_x s_k&=0.&\label{entropy subsystems: convective}
    \end{alignat}
\end{subequations}
Basically, the splitting decouples the Lagrangian derivative terms from the remaining terms. Therefore, the convective system can be written as ${\rm D} Q/{\rm D} t =0$ for $Q\in\left\{\rho, \rho u, \rho E, \rho Y_1, \alpha_1\right\}$. Now, the acoustic system contains all the pressure terms and the interfacial term of the topological equation (\ref{model volume fraction eq}). Note that this interfacial term includes the spatial derivative of velocity and is therefore included in the acoustic system. \added[id=Rev.1]{The splitting step is first-order accurate in time. A higher-order temporal accuracy can be obtained, e.g. for second-order accuracy by employing Strang splitting \cite{leveque2002finite}.}\\
\indent The numerical solution of (\ref{governing equations 5 equation model}) consists of successively approximating the solution of the acoustic system and the convective system. By denoting the temporal step size with $\Delta t$, the mesh width with $\Delta x$, the fluid state at time $n \Delta t$ and position $j \Delta x$ with $\mathbf{Q}_j^n\equiv(\rho, \rho u, \rho E, \rho Y_1, \alpha_1)_j^n$, and an intermediate time level with $n+1-$, the approximation within one time step reads:
\begin{enumerate}
 \item Update $\mathbf{Q}_j^n$ to $\mathbf{Q}_j^{n+1-}$ by time marching the acoustic system (\ref{acoustic}) with step size $\Delta t$;
 \item Update $\mathbf{Q}_j^{n+1-}$ to $\mathbf{Q}_j^{n+1}$ by time marching the convective system (\ref{convective}) with step size $\Delta t$.
\end{enumerate}
The choice of numerically solving the submodels in this order is linked to the velocity approximation: the velocity of the acoustic system is used for the determination of the convective velocity in order to ensure the conservation of mass, momentum, energy and partial masses as it is detailed in Section \ref{Main properties scheme}. The details of each step are given in Section\added[id=Authors]{s} \ref{Acoustic system} and \ref{convective system}.
\subsection{Mathematical analysis of the two submodels}\label{subsection Mathematical analysis}


The five\added[id=Authors]{-}equation model (\ref{governing equations 5 equation model}) can be cast into the primitive form
    \begin{equation}
        \partial_t\mathbf{W}+\mathbf{B(W)}\partial_x\mathbf{W}=\mathbf{0},
    \end{equation}
and the primitive form of the subsystems (\ref{acoustic})-(\ref{convective}) reads: (i) for the acoustic system:
    \begin{equation}
        \partial_t\mathbf{W}+\mathbf{A(W)}\partial_x\mathbf{W}=\mathbf{0},
    \end{equation}
and (ii) for the convective system:
    \begin{equation}
        \partial_t\mathbf{W}+\mathbf{C(W)}\partial_x\mathbf{W}=\mathbf{0},
    \end{equation}
where
    \begin{equation}
        \begin{array}{l}
            \mathbf{B(W)}=\mathbf{A(W)}+\mathbf{C(W)},
        \end{array}
    \end{equation}
with
\begin{equation}
\mathbf{W}=\begin{pmatrix}\rho\\[6pt]u\\[6pt]p\\[6pt] Y_1\\[6pt] \alpha_1\end{pmatrix},~~\mathbf{A(W)}=\begin{pmatrix}0 & \rho & 0 & 0 & 0\\[6pt]0 & 0 & 1/\rho & 0 & 0\\[6pt]0 & \rho c^2 & 0 & 0 & 0\\[6pt]0 & 0 & 0 & 0 & 0 \\[6pt] 0 & K & 0 & 0 & 0\end{pmatrix}, \quad \mathbf{C(W)}=u~\mathbf{I}_5,\end{equation}
where $\mathbf{I}_d$ is the identity matrix in $\mathbb{R}^{d\times d}$. The derivation of the pressure equation is straightforward and can be found in \cite{murrone2005five, kreeft2010new, ten2015compressible}. This casting reveals that the matrix $\mathbf{B}$ splits into an acoustic part $\mathbf{A}$ and a convective part $\mathbf{C}$. The eigenvalues of the full system ($\lambda_k$) split also into an acoustic part ($\lambda_k^a$) and a convective part ($\lambda_k^c$) as $\lambda_k=\lambda_{k}^a+\lambda_{k}^c$ with
    \begin{equation}\label{Acoustic system eigenvalues}
        \begin{array}{l l l}
            \lambda_{1}^a=-c,& \lambda_{2,3,4}^a=0,& \lambda_{5}^a=c,\\[8pt]
            \lambda_{1}^c=u,& \lambda_{2,3,4}^c=u,& \lambda_{5}^c=u.
        \end{array}
    \end{equation}
The characteristic fields associated with the convective submodel are obviously LD. Concerning the acoustic submodel, the fields associated with the middle wave $\lambda^{a}_{2,3,4}=0$ are LD. The other two waves, associated with $\lambda^{a}_1=-c, \lambda^{a}_5=c$, can be shown, by using a similar argument as Murrone et al. \cite{murrone2005five}, to be GNL in the non-isobaric case and LD in the isobaric case.


\subsection{Numerical solution of the acoustic submodel}\label{Acoustic system}


\subsubsection{Lagrangian formulation}\label{Acoustic system: New formulation}


Introducing the specific volume $\tau=1/\rho$ and taking $\left\{\tau, u, E, Y_1, \alpha_1\right\}$ as the set of variables, the acoustic system can be cast into the form
    \begin{subequations}\label{acoustic HLLC}
        \begin{alignat}{6}
            &\partial_t\tau&-&\tau \partial_xu&=0,\label{1acoustic rho}\\
            &\partial_tu&+&\tau  \partial_xp&=0,\label{1acoustic rho u}\\
            &\partial_tE&+&\tau\partial_x(pu)&=0,\label{1acoustic rho E}\\
            &\partial_tY_1 &&&=0,\label{1acoustic alpha1 rho1}\\
            &\partial_t\alpha_1&+&\rho K \tau\partial_xu&=0.\label{1acoustic alpha1 rho1 E1}
        \end{alignat}
    \end{subequations}
The Eqs. (\ref{1acoustic rho})-(\ref{1acoustic rho E}) describe the bulk fluid, and the Eqs. (\ref{1acoustic alpha1 rho1})-(\ref{1acoustic alpha1 rho1 E1}) describe the evolution of the fraction variables, which are specific for the five-equation two-phase flow model. The second term of each equation (except the fourth) contains the operator $\tau \partial_x$. \added[id=Rev.1]{As in \cite{chalons2014all}, for $t \in [t^n,t^n+\Delta t)$ we approximate $\tau(x,t)\partial_x$ by $\tau(x,t^n)\partial_x$, where the time level is $t^n=n\Delta t$ with time step $\Delta t$. We then introduce the mass variable $m$ by ${\rm d} m =\rho(x,t^n) {\rm d}x $. The \textit{Lagrangian} system}
    \begin{subequations}\label{acoustic HLLC2}
        \begin{alignat}{6}
            &\partial_t\tau&-& \partial_mu&=0,\\
            &\partial_tu&+&  \partial_mp&=0,\\
            &\partial_tE&+&\partial_m(pu)&=0,\\
            &\partial_tY_1 &&&=0,\\
            &\partial_t\alpha_1&+&\rho K \partial_mu&=0,
        \end{alignat}
    \end{subequations}
is a first-order in time approximation of (\ref{acoustic HLLC}). This system has the eigenvalues
    \begin{equation}\label{Acoustic system eigenvalues lag}
        \begin{array}{l l l}
            (\lambda_{1}^a)^{\mathcal{L}ag}=-\rho c,& (\lambda_{2,3,4}^a)^{\mathcal{L}ag}=0,& (\lambda_{5}^a)^{\mathcal{L}ag}=\rho c
        \end{array}
    \end{equation}
\added[id=Rev.2]{and associated eigenvectors}
    \begin{equation}\label{Acoustic system eigenvectors}
        \begin{array}{l l l}
            (\mathbf{v}_{1}^a)^{\mathcal{L}ag}=\begin{pmatrix}  -1\\  -\rho c\\  (\rho c)^2\\  0\\  \rho K\end{pmatrix},\quad  (\mathbf{v}_{2}^a)^{\mathcal{L}ag}=\begin{pmatrix} 1\\ 0\\ 0\\ 0\\ 0\end{pmatrix}, \quad (\mathbf{v}_{3}^a)^{\mathcal{L}ag}=\begin{pmatrix} 0\\ 0\\ 0\\ 1\\ 0\end{pmatrix},\\[35pt]
             (\mathbf{v}_{4}^a)^{\mathcal{L}ag}=\begin{pmatrix} 0\\ 0\\ 0\\ 0\\ 1\end{pmatrix}, \quad (\mathbf{v}_{5}^a)^{\mathcal{L}ag}=\begin{pmatrix} -1\\  \rho c\\  (\rho c)^2\\  0\\  \rho K\end{pmatrix}.
        \end{array}
    \end{equation}
It can be written in the following vectorial form:
    \begin{equation}\label{a 5}
        \partial_t\mathbf{Q}^{\mathcal{L}ag}+\partial_m\mathcal{F}^{\mathcal{L}ag}(\mathbf{Q}^{\mathcal{L}ag})+\mathcal{B}^{\added[id=Authors]{\mathcal{L}ag}}\left(\mathbf{Q}^{\mathcal{L}ag}\right)\partial_m u=\mathbf{0},
    \end{equation}
where
    \begin{subequations}\label{relaxed system 5}
        \begin{alignat}{2}
            \mathbf{Q}^{\mathcal{L}ag}&=(\tau, u, E, Y_1, \alpha_1)^T,\\
            \mathcal{F}^{\mathcal{L}ag}(\mathbf{Q}^{\mathcal{L}ag})&=(-u, p, p u, 0, 0)^T,\\
            \mathcal{B}^{\mathcal{L}ag}(\mathbf{Q}^{\mathcal{L}ag})&=(0, 0, 0, 0,\rho K)^T.\label{Blag}
        \end{alignat}
    \end{subequations}
The superscript $\mathcal{L}ag$ is used for the variables in the Lagrangian system. The term $\mathcal{F}^{\mathcal{L}ag}$ is a conservative flux and the latter is the non-conservative term. \added[id=Rev.2]{System (\ref{a 5})-(\ref{relaxed system 5}) is numerically approximated in the following.}


\subsubsection{HLLC-type solver for the acoustic submodel in Lagrangian coordinates}\label{subsubsection HLLC-type solver}

\begin{figure}[h!]
\centering
\includegraphics[width=0.4\textwidth]{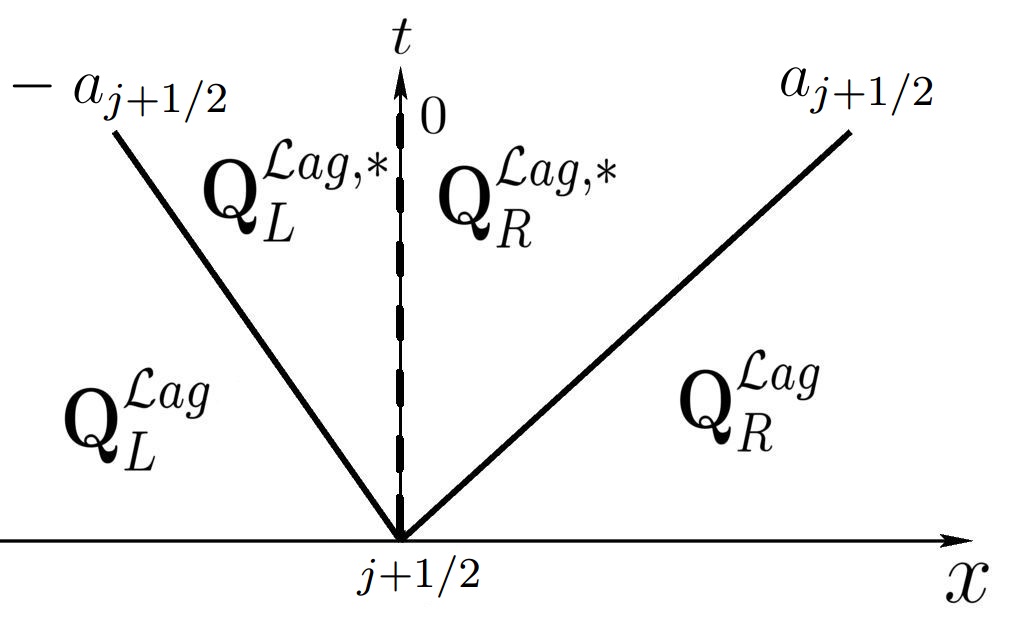}
\caption{The different states $\mathbf{Q}^{\mathcal{L}ag}_L, \mathbf{Q}^{\mathcal{L}ag,*}_L, \mathbf{Q}^{\mathcal{L}ag,*}_R, \mathbf{Q}^{\mathcal{L}ag}_R$ and wave speeds $-a_{j+1/2}, 0, a_{j+1/2}$ in the Riemann problem.}
\label{Riemann_problem}
\end{figure}

An HLLC-type Riemann solver \cite{toro1994restoration} is used to solve the acoustic system (\ref{a 5})-(\ref{relaxed system 5}). The finite-volume approximation of the Eqs. (\ref{a 5})-(\ref{relaxed system 5}) on each mesh element $[x_{j-1/2},x_{j+1/2}]$ \added[id=Rev.1]{ follows from integration over the mesh element and assuming a constant density in the $m$ variable and constant interfacial term in each element}, and reads
\begin{equation}\label{rightmost}
    \begin{array}{l l}
        \partial_t(\left(\mathbf{Q}^{\mathcal{L}ag}\right)_j)&+\dfrac{1}{\Delta m_j}\left(\left(\boldsymbol{F}^{\mathcal{L}ag}\right)^{\text{HLLC}}_{j+1/2}-\left(\boldsymbol{F}^{\mathcal{L}ag}\right)^{\text{HLLC}}_{j-1/2}\right)
        \\[10pt]&+\added[id=Rev.1]{\dfrac{1}{\Delta m_j}}\mathcal{B}^{\added[id=Authors]{\mathcal{L}ag}}\added[id=Rev.1]{\left(\left(\mathbf{Q}^{\mathcal{L}ag}\right)_j\right)\left(u^*_{j+1/2}-u^*_{j+1/2}\right)}=0,
    \end{array}
\end{equation}
\added[id=Rev.1]{with $\Delta m_j=\rho_j^n\Delta x$.} \added[id=Rev.2]{In this paper we employ the classical finite-volume notation in which subscript $j$ refers to a cell average and $j+1/2$ to a cell boundary.} The HLLC-type numerical flux vector $\left(\boldsymbol{F}^{\mathcal{L}ag}\right)^{\text{HLLC}}$, which approximates $\mathcal{F}^{\mathcal{L}ag}\left(\mathbf{Q}^{\mathcal{L}ag}\right)$, is obtained by applying the HLLC-type relations across the three different waves with eigenvalues (\ref{Acoustic system eigenvalues lag}), see Figure \ref{Riemann_problem}.
\added[id=Rev.2]{Using (\ref{Acoustic system eigenvectors}) we see that the velocity and pressure are the Riemann invariants of the LD middle wave}. The HLLC-type relations across the left and right waves for the momentum equation are given by
    \begin{subequations}
        \begin{alignat}{2}
            p^*_{j+1/2}&=p_{j}-a_{j+1/2} \left(u_{j+1/2}^*-u_j\right),\\
            p^*_{j+1/2}&=p_{j+1}+a_{j+1/2} \left(u_{j+1/2}^*-u_{j+1}\right).
        \end{alignat}
    \end{subequations}
where the acoustic impedance $a_{j+1/2}$ at the interface is estimated using the direct computation of the eigenvalues of the acoustic submodel\added[id=Authors]{:}
    \begin{equation}
        a_{j+1/2}=\max (\rho_j c_j, \rho_{j+1} c_{j+1}).
    \end{equation}
This leads to a single-state HLLC numerical flux-vector:
    \begin{equation}
        \left(\boldsymbol{F}^{\mathcal{L}ag}\right)^{\text{HLLC}}_{j+1/2}= \left(-u^*, p^*, p^* u^*, 0, 0\right)_{j+1/2},
    \end{equation}
where
    \begin{subequations}
        \begin{alignat}{2}
            u^*_{j+1/2}&=\dfrac{u_j+u_{j+1}}{2}+\dfrac{p_j-p_{j+1}}{2a_{j+1/2}},\\
            p^*_{j+1/2}&=\dfrac{p_j+p_{j+1}}{2}+\frac{a_{j+1/2}}{2}(u_j-u_{j+1}).
        \end{alignat}
    \end{subequations}
The interfacial term of the topology equation is approximated at first-order by
    \begin{equation}\label{approximation integral u 1}
        \added[id=Authors]{K^n_j \left(u_{j+1/2}^*-u_{j-1/2}^*\right).}
    \end{equation}
Summarizing and using an explicit forward Euler time step, the update formula for the discretized acoustic system reads:
    \begin{equation}\label{Update formulae4}
        \begin{array}{l l}
            \left(\mathbf{Q}^{\mathcal{L}ag}\right)_j^{n+1-}=\left(\mathbf{Q}^{\mathcal{L}ag}\right)_j^{n}&-\dfrac{\Delta t}{\rho_j^n\Delta x}\left(\left(\boldsymbol{F}^{\mathcal{L}ag}\right)^{\text{HLLC},n}_{j+1/2}-\left(\boldsymbol{F}^{\mathcal{L}ag}\right)^{\text{HLLC},n}_{j-1/2}\right)\\
            &-K_j^n\dfrac{\Delta t}{\Delta x}\left(\left(\mathcal{H}^{\mathcal{L}ag}\right)^{n}_{j+1/2}-\left(\mathcal{H}^{\mathcal{L}ag}\right)^{n}_{j-1/2}\right),
        \end{array}
    \end{equation}
where
    \begin{equation}\label{Update formulae: H}
            \left(\mathcal{H}^{\mathcal{L}ag}\right)^{T}=(0, 0, 0, 0, u^*).
    \end{equation}
\added[id=Rev.2]{The numerical experiments in section \ref{section Numerical results} employ this update formula.}

\subsubsection{Update of the acoustic submodel in Eulerian variables}


The update formulae for the discretized acoustic system in terms of the Eulerian variables from (\ref{governing equations 5 equation model}) are \added[id=Rev.1]{a reformulation of those in (\ref{Update formulae4})-(\ref{Update formulae: H}) and read}:
    \begin{subequations}\label{update formulae acoustic system 1}
        \begin{alignat}{2}
            R_j\rho_j^{n+1-}&=\rho_j^n,\\
            R_j\left(\rho u\right)_j^{n+1-}&=\left(\rho u\right)_j^n-\frac{\Delta t}{\Delta x}\left(p^*_{j+1/2}-p^*_{j-1/2}\right),\\
            R_j\left(\rho E\right)_j^{n+1-}&=\left(\rho E\right)_j^n-\frac{\Delta t}{\Delta x}\left(p^*_{j+1/2}u^*_{j+1/2}-p^*_{j-1/2}u^*_{j-1/2}\right),\\
            R_j\left(\rho Y_1\right)_j^{n+1-}&=\left(\rho Y_1\right)_j^n,\label{update acoustic Y1 mass fraction}\\
            \left(\alpha_1\right)_j^{n+1-}&=\left(\alpha_1\right)_j^n-K_j^n\frac{\Delta t}{\Delta x}\left(u^*_{j+1/2}-u^*_{j-1/2}\right),\label{update acoustic alpha1}
        \end{alignat}
    \end{subequations}
where $R_j$ is given by
    \begin{equation}\label{def Rj}
        R_j=1+\frac{\Delta t}{\Delta x}\left(u^*_{j+1/2}-u^*_{j-1/2}\right).
    \end{equation}
\added[id=Rev.2]{Some properties of the numerical scheme, presented in section \ref{Main properties scheme}, employ these update formulae in the derivation.}
\subsection{Numerical solution of the convective submodel}\label{convective system}
The convective system is approximated by using a classical upwind finite-volume scheme as employed in Chalons et al. \cite{chalons2014all}. Making again a forward Euler time step, the scheme reads:
    \begin{equation}\label{alternative differential form: eq energy fluid 1}
        \begin{array}{l l}
            \varphi_j^{n+1}=\varphi_j^{n+1-}&-\dfrac{\Delta t}{\Delta x}\left(u^*_{j+1/2}\varphi_{j+1/2}^{n+1-}-u^*_{j-1/2}\varphi_{j-1/2}^{n+1-}\right)\\[10pt]
            &+\dfrac{\Delta t}{\Delta x}\varphi_{j}^{n+1-}\left(u^*_{j+1/2}-u^*_{j-1/2}\right).
        \end{array}
    \end{equation}
where $\varphi\in\left\{\rho, \rho u, \rho E, \rho Y_1, \alpha_1\right\}$. The upwind value is used to approximate the interface value $\varphi_{j+1/2}$:
    \begin{equation}
        \varphi_{j+1/2}^{n+1-}=
            \left\{
                \begin{array}{l l}
                    \varphi_{j}^{n+1-},  & \quad \text{if}\quad u_{j+1/2}^*\geq 0,\\[8pt]
                    \varphi_{j+1}^{n+1-},  & \quad \text{if}\quad u_{j+1/2}^*<0.
                \end{array}
            \right.
    \end{equation}


\subsection{Stability requirement}


The common time step in the explicit time integration method is obtained using the Courant numbers of both subsystems. The Courant numbers are given by
    \begin{equation}\label{stability acoustic system}
        \mathcal{C}^{a}=\dfrac{\Delta t}{\Delta x}\max_j \lambda_j^a,
    \end{equation}
with maximum wave speed $\lambda_j^a=\max\left(\tau_j^n,\tau_{j+1}^n\right)a_{j+1/2}$, for the acoustic subsystem, and by
    \begin{equation}\label{CFL convective system}
       \mathcal{C}^{c}=\dfrac{\Delta t}{\Delta x} \max_j \lambda_j^c,
    \end{equation}
with the maximum wave speed $\lambda_j^c=\left(u_{j-1/2}^*\right)^+-\left(u_{j+1/2}^*\right)^-$, for the convective subsystem, where $b^{\pm}=(b\pm |b|)/2$. The time step is determined by the requirement that both Courant numbers need to be less than one. In the implementation, the most severe time step restriction is taken for both subsystems. Hence, the time step size is selected with the Courant number $\mathcal{C}=\max\left(\mathcal{C}^a, \mathcal{C}^c\right)$. The Courant number for the classical direct approaches is defined by
    \begin{equation}\label{CFL convective system}
       \mathcal{C}^{d}=\dfrac{\Delta t}{\Delta x} \max_j \left(|u_{j+1/2}|+c_{j+1/2}\right).
    \end{equation}


\subsection{\added[id=Rev.3]{Main properties scheme}}\label{Main properties scheme}

\subsubsection{Conservation of mass, momentum, energy and partial mass}
The scheme of the convective system (\ref{alternative differential form: eq energy fluid 1}) can be written as:
    \begin{equation}\label{alternative differential form: eq energy fluid 12}
         \varphi_j^{n+1}=R_j\varphi_j^{n+1-}-\dfrac{\Delta t}{\Delta x}\left(u^*_{j+1/2}\varphi_{j+1/2}^{n+1-}-u^*_{j-1/2}\varphi_{j-1/2}^{n+1-}\right),
    \end{equation}
where $R_j$ is defined by (\ref{def Rj}). Substitution of (\ref{update formulae acoustic system 1}) into this form leads to the update formulae
    \begin{subequations}\label{preservation of mass, momentum, energy and partial mass}
        \begin{alignat}{2}
            (\rho)_j^{n+1}      =&(\rho )_j^{n}     -\frac{\Delta t}{\Delta x}\left(u^*_{j+1/2}\rho_{j+1/2}^{n+1-}-u^*_{j-1/2}\rho_{j-1/2}^{n+1-}\right),\\
            (\rho u)_j^{n+1}    =&(\rho u)_j^{n}    \nonumber \\
                                 &-\frac{\Delta t}{\Delta x}\left(u^*_{j+1/2}(\rho u)_{j+1/2}^{n+1-}+p^*_{j+1/2}-u^*_{j-1/2}(\rho u)_{j-1/2}^{n+1-}-p^*_{j-1/2}\right),\\
            (\rho E)_j^{n+1}   =&(\rho E)_j^{n}    -\frac{\Delta t}{\Delta x}\left(u^*_{j+1/2}(\rho E)_{j+1/2}^{n+1-}+p^*_{j+1/2}u^*_{j+1/2}\right.\nonumber \\
                                 &\left. \quad \quad -u^*_{j-1/2}(\rho E)_{j-1/2}^{n+1-}-p^*_{j-1/2}u^*_{j-1/2}\right),\\
            (\rho Y_1)_j^{n+1}  =&(\rho Y_1)_j^{n}  -\frac{\Delta t}{\Delta x}\left(u^*_{j+1/2}(\rho Y_1)_{j+1/2}^{n+1-}-u^*_{j-1/2}(\rho Y_1)_{j-1/2}^{n+1-}\right),
        \end{alignat}
    \end{subequations}
which guarantees the conservation of mass, momentum, energy and partial mass of the proposed approach. Please notice that the choice of $u^*_{j+1/2}$ in the transport scheme makes it possible to have a fully conservative scheme for the conservative variables \cite{chalons2014all}. Due to the non-conservative form of the topology equation, there is no conservation of the volume fraction.
\subsubsection{Positivity of the volume fraction and mass fraction}
Using the definition of the interfacial variable (\ref{K}), the update formula (\ref{update acoustic alpha1}) of the volume fraction in the acoustic system can be written as
\begin{equation}\label{update alpha rewritten}
\begin{array}{l l}
    \left(\alpha_1\right)_j^{n+1-}&=\left(\alpha_1\right)_j^n\left[1-\dfrac{\Delta t}{\Delta x}\left(\alpha_2\right)_j^n\dfrac{\left(\rho_2c_2^2\right)_j^n-\left(\rho_1 c_1^2\right)_j^n}{\left(\alpha_2\right)_j^n\left(\rho_1c_1^2\right)_j^n+\left(\alpha_1\right)_j^n\left(\rho_2c_2^2\right)_j^n}\left(u^*_{j-1/2}-u^*_{j+1/2}\right)\right].
    \end{array}
\end{equation}
Since $\left(\alpha_1\right)_j^n \geq 0$, positivity of the volume fraction is ensured when the part within the brackets is positive, i.e.
\begin{equation}\label{positivity condition}
    A_j^n\frac{\Delta t}{\Delta x}\left(u^*_{j-1/2}-u^*_{j+1/2}\right)\leq 1,
\end{equation}
where \begin{equation}\label{def A}
    A_j^n=\left(\alpha_2\right)_j^n\dfrac{\left(\rho_2c_2^2\right)_j^n-\left(\rho_1 c_1^2\right)_j^n}{\left(\alpha_2\right)_j^n\left(\rho_1c_1^2\right)_j^n+\left(\alpha_1\right)_j^n\left(\rho_2c_2^2\right)_j^n}.
\end{equation}
The observations
\begin{subequations}\label{update formulae acoustic system 1 conservative}
        \begin{alignat}{1}
            \left(\rho_2 c_2^2\right)_j^n-\left(\rho_1c_1^2\right)_j^n<\max\left[\left(\rho_1 c_1^2\right)_j^n,\left(\rho_2c_2^2\right)_j^n\right],\\
            \left(\alpha_2\right)_j^n\left(\rho_1 c_1^2\right)_j^n+\left(\alpha_1\right)_j^n\left(\rho_2c_2^2\right)_j^n>\min\left[\left(\rho_1 c_1^2\right)_j^n,\left(\rho_2c_2^2\right)_j^n\right],
        \end{alignat}
    \end{subequations}
and $0 \leq \left(\alpha_2\right)_j^n \leq 1$ imply that $A_j^n \leq 1$. Using the CFL-type condition given in (\ref{CFL convective system}), we obtain 
\begin{equation}\label{positivity condition CFL}
    \frac{\Delta t}{\Delta x}\left(u^*_{j-1/2}-u^*_{j+1/2}\right)\leq \frac{\Delta t}{\Delta x}\left[\left(u^*_{j-1/2}\right)^+-\left(u^*_{j+1/2}\right)^-\right]\leq 1.
\end{equation}
Positivity of the volume fraction is thus ensured by combining the results. Note that the upper bound $\left(\alpha_1\right)_j^{n+1-} \leq 1$ is a direct consequence of this result. Similarly, the update formula (\ref{update acoustic Y1 mass fraction}) ensures the positivity of the mass fraction.
%
%

\section{Numerical results}\label{section Numerical results}


\added[id=Rev.2]{To illustrate the behavior of the proposed scheme, it is evaluated for} \added[id=Rev.3]{five two-phase flow} problems encountered in the literature: a translating interface problem, a pressure jump problem, a no-reflection problem\added[id=Authors]{, a water-air mixture} problem and a \added[id=Rev.3]{two-phase cavitation problem}. To illustrate the behavior of the proposed scheme, we consider standard shock-tube problems encountered in the literature.\\
All \added[id=Rev.3]{five} test cases are defined such that no wave hits a boundary before the prescribed end time. All test cases are also computed using the direct HLLC-type approach proposed by Daude et al. \cite{daude2014numerical}. The tests are performed with first\added[id=Authors]{-}order accuracy in space and time. For each test, the Courant numbers of the current splitting approach and the direct approach are taken equal: $\mathcal{C}=\mathcal{C}^d$. The comparisons are performed using the same number of cells. \added[id=Revs]{The convergence rates are shown for each test case where an analytical solution is available.} \added[id=Rev.1]{To compare the performance of both methods, the number of time steps and the CPU times are reported.}


\subsection{Translating two-phase interface}


In this first test case, also considered in e.g. \cite{kreeft2010new}, a dense fluid and a much less dense gas move to the right, at constant velocity and pressure. The initial interface is located in the middle of the tube ($x=0.0$) of length $L=0.5$. This test case is considered to assess the behavior of the present scheme at a material interface with a density jump which is representative for that of the important class of water-air flows.\\
\begin{table}[h!]
\caption{Initial values and material properties for the translating interface problem.}\label{init translating interface}
\begin{tabu} to \textwidth {X X X X X X X X X X}
   \toprule
   \multicolumn{3}{c}{\footnotesize{(a) Initial values}} & &&&& \multicolumn{2}{l}{\footnotesize{(b) Material properties}}\\
   \cline{1-6} \cline{8-9}
   & \footnotesize{$ \rho $}   & \footnotesize{$ u $}     & \footnotesize{$ p $}   & \footnotesize{$ Y_1 $} & \footnotesize{$\alpha_1$} & & & \footnotesize{$\gamma$  } \\\hline
   \footnotesize{Fluid 1}   &   \footnotesize{$ 1000 $}   & \footnotesize{$ 1.0 $  } & \footnotesize{$ 1.0 $} & \footnotesize{$ 1.0  $} & \footnotesize{$1.0$  } &&\footnotesize{Fluid 1}  &  \footnotesize{$1.4$  } \\
   \footnotesize{Fluid 2}   &   \footnotesize{$ 1.0 $}   & \footnotesize{$ 1.0 $  } & \footnotesize{$ 1.0 $} & \footnotesize{$ 0.0  $} & \footnotesize{$0.0$  } &&\footnotesize{Fluid 2}  &  \footnotesize{$1.6$  } \\
   \bottomrule
\end{tabu}
\end{table}\\

\begin{figure}[h!]
\centering
\includegraphics[width=0.6\textwidth]{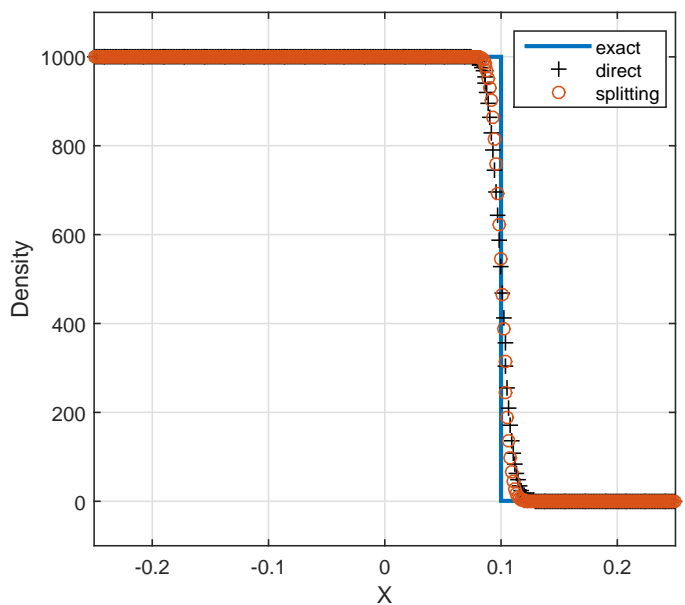}
\caption{Translating interface problem - density profile - Exact solution ``{\color{blue}-}'', splitting approach ``{\color{Orange}$\boldsymbol{\circ}$}'' and direct approach ``$\bold{+}$"at $t=0.1$.}
\label{Translating_interface_nr2_density}
\end{figure}

\begin{figure}[h!]
\centering
\includegraphics[width=0.6\textwidth]{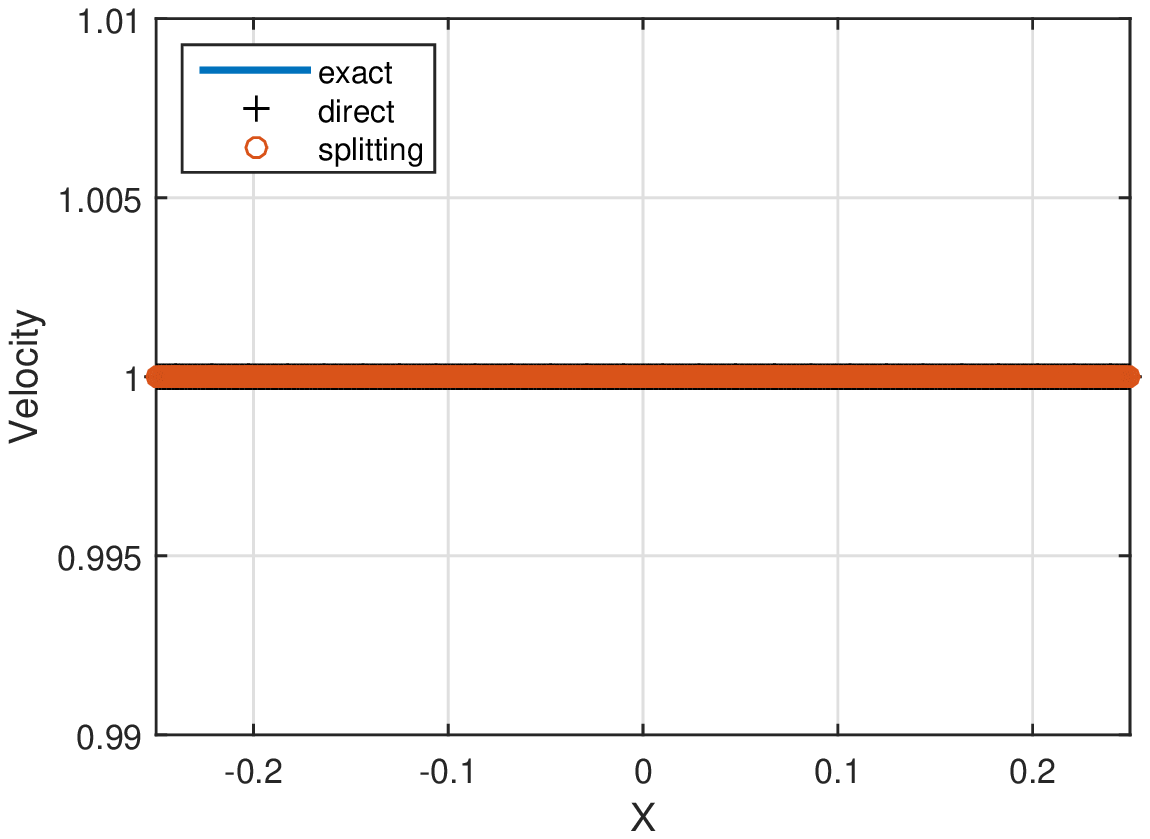}
\caption{Translating interface problem - velocity profile - Exact solution ``{\color{blue}-}'', splitting approach ``{\color{Orange}$\boldsymbol{\circ}$}'' and direct approach ``$\bold{+}$"at $t=0.1$.}
\label{Translating_interface_nr2_velocity}
\end{figure}

\begin{figure}[h!]
\centering
\includegraphics[width=0.6\textwidth]{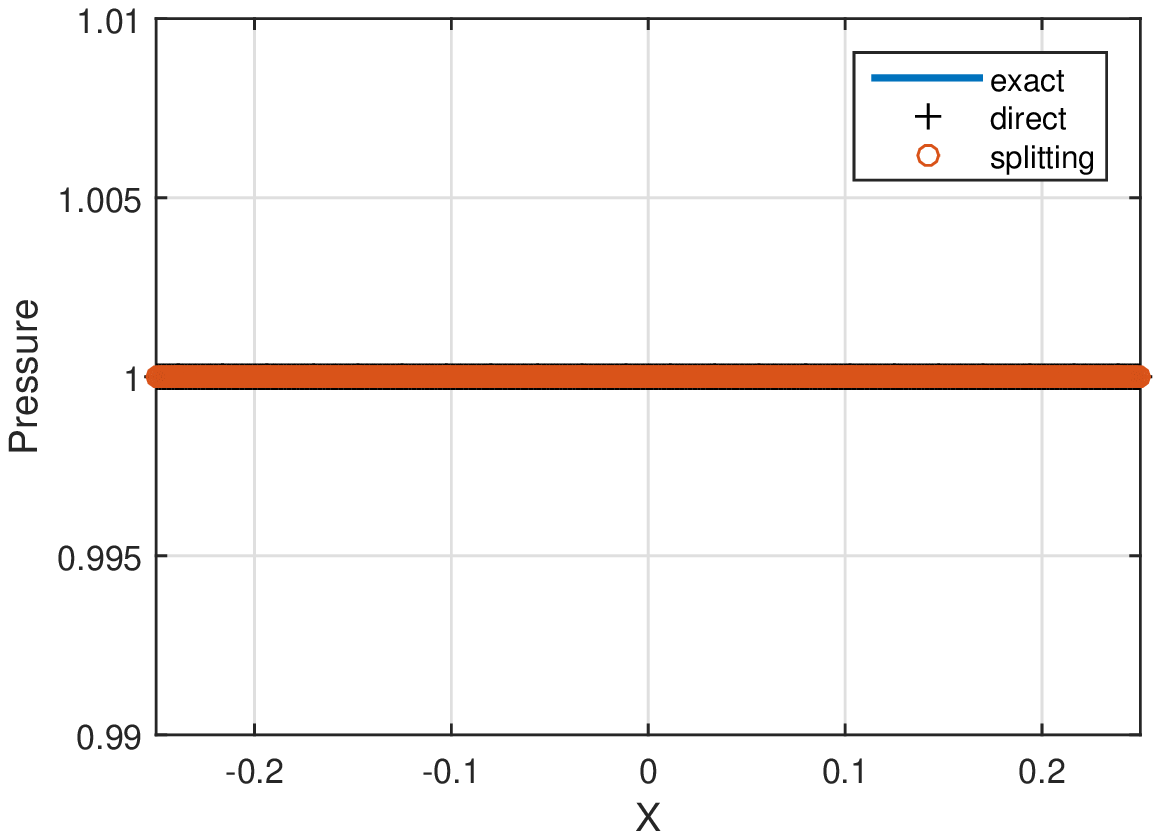}
\caption{Translating interface problem - pressure profile  - Exact solution ``{\color{blue}-}'', splitting approach ``{\color{Orange}$\boldsymbol{\circ}$}'' and direct approach ``$\bold{+}$"at $t=0.1$.}
\label{Translating_interface_nr2_pressure}
\end{figure}

\begin{figure}[h!]
\centering
\includegraphics[width=0.6\textwidth]{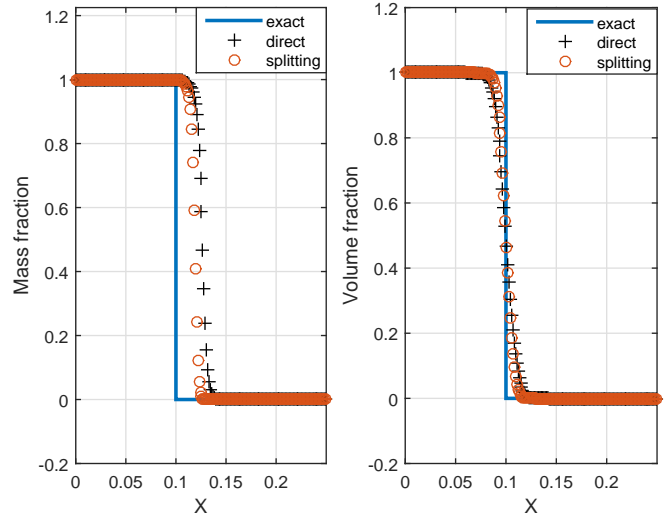}
\caption{Translating interface problem - volume and mass fraction profiles  - Exact solution ``{\color{blue}-}'', splitting approach ``{\color{Orange}$\boldsymbol{\circ}$}'' and direct approach ``$\bold{+}$"at $t=0.1$.}
\label{Translating_interface_nr2_fractions}
\end{figure}

The initial values and material properties are given in Table \ref{init translating interface}. Two perfect gases are considered ($\pi_1=\pi_2=0$, \added[id=Authors]{$\eta_1=\eta_2=0$}), with the difference for both fluids only in $\gamma$. The depicted results have been obtained at time $t=0.1$ with $N=400$ cells and a Courant number $\mathcal{C}=0.95$. The distributions of the primitive variables are visualized in the Figures \ref{Translating_interface_nr2_density} to \ref{Translating_interface_nr2_fractions} \added[id=Revs]{and the convergence rates of the density profiles are listed in Table \ref{convratesnr2}.}
\begin{figure}[h!]
\centering
\includegraphics[width=1.0\textwidth]{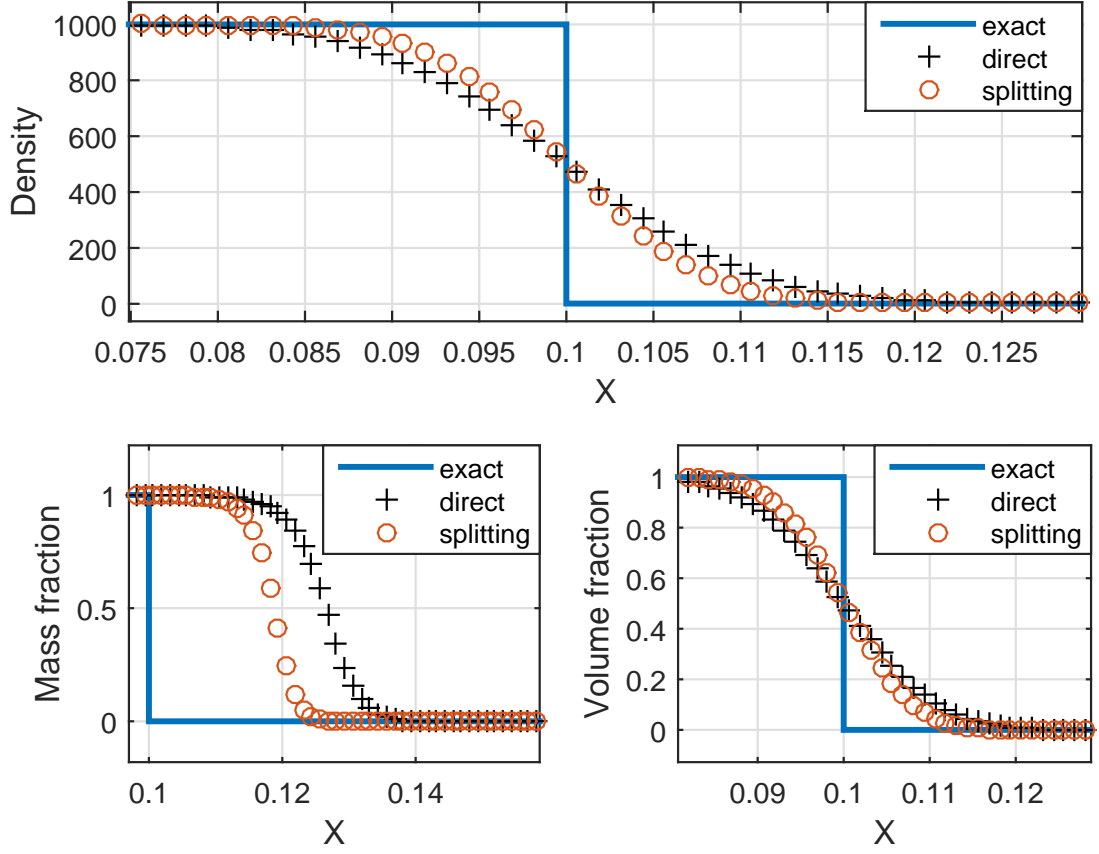}
\caption{Translating interface problem - zoom at contact discontinuity - Exact solution ``{\color{blue}-}'', splitting approach ``{\color{Orange}$\boldsymbol{\circ}$}'' and direct approach ``$\bold{+}$"at $t=0.1$.}
\label{Translating_interface_nr2_zoom}
\end{figure}\\

\begin{table}[h!]
\caption{The $L_1$-convergence rates for the density of the translating interface problem. The convergence rates are computed as $c_N=\log(e_N/e_{2N})/\log(2)$. The errors are given by $e_N=\|s_N-s_{\text{exact}}\|_{L_1}$, where $s_N$ is the solution computed with $N$ grid points, $s_{\text{exact}}$ the exact solution, and $\|\cdot\|_{L_1}$ the standard $L_1$-norm.}\label{convratesnr2}
\begin{tabu} to \textwidth {X X X }
   \toprule
   \footnotesize{Convergence rates}& \footnotesize{Splitting}   & \footnotesize{Direct}   \\ \hline
   \footnotesize{$c_{40}$}   &   \footnotesize{$ 0.67 $}   & \footnotesize{$ 0.56 $  }  \\
   \footnotesize{$c_{80}$}   &   \footnotesize{$ 0.64 $}   & \footnotesize{$ 0.53 $  }   \\
   \footnotesize{$c_{160}$}  &   \footnotesize{$ 0.63 $}   & \footnotesize{$ 0.52 $  }  \\
   \footnotesize{$c_{320}$}  &   \footnotesize{$ 0.60 $}   & \footnotesize{$ 0.51 $  }  \\
   \footnotesize{$c_{640}$}  &   \footnotesize{$ 0.57 $}   & \footnotesize{$ 0.50 $  }  \\
      \bottomrule
\end{tabu}
\end{table}

The results obtained with the proposed splitting-based method are very similar to the ones obtained with the direct approach from \cite{daude2014numerical}. The contact discontinuity is well retrieved with both methods, whereas the velocity and pressure profiles are perfectly constant; no pressure oscillations occur across the interface. The location of \added[id=Authors]{the} two-phase interface for the mass fraction is a bit off (see Figure \ref{Translating_interface_nr2_zoom}), for both the proposed method and the direct approach from \cite{daude2014numerical}. This is also the case for the method proposed in \cite{kreeft2010new}. In the region where the material interface is smeared due to intrinsic numerical dissipation of the two numerical schemes, the associated cells contain both fluids with $\alpha_2\rho_2 \ll \alpha_1 \rho_1$ which gives a value of $Y_1$ close to $1$. \added[id=Authors]{With a finer mesh, the correct location is obtained, see also Table \ref{convratesnr2}.} At the end time $t=0.1$ the contact discontinuity is indeed located at $x=0.1$. \added[id=Authors]{The proposed method captures the location slightly better.}  The newly proposed method takes larger time steps ($124$ time steps) than the direct approach from \cite{daude2014numerical} ($192$ time steps). \added[id=Rev.1]{The CPU time is $0.17$s and $0.36$s for the splitting approach and the direct approach, respectively (averaged over $500$ runs on an i5 processor).} \added[id=Authors]{Both methods show similar convergence rates, see Table \ref{convratesnr2}.}
\clearpage

\subsection{A two-pressure jump problem}


In this test case, proposed by Barberon al. \cite{barberon2003practical} and also considered in \cite{daude2014numerical}, the shock tube is \added[id=Authors]{again} filled with two \added[id=Authors]{perfect} gases with different densities. The pressures at both sides are slightly different. The interface is located at $x=0.5$ m. Due to the pressure difference, a shock wave will propagate rightwards and a rarefaction wave will propagate leftwards.\\

\begin{table}[h!]
\caption{Initial values and material properties for the two-pressure jump problem. The dimensions of the quantities $\rho, u$ and $p$ are kg m$^{-3}$, m s$^{-1}$ and Pa, respectively.}\label{init pressure jump}
\begin{tabu} to \textwidth {X X X X X X X X X X}
   \toprule
   \multicolumn{3}{c}{\footnotesize{(a) Initial values}} & &&&& \multicolumn{2}{c}{\footnotesize{(b) Material properties}}\\
   \cline{1-6} \cline{8-9}
   & \footnotesize{$ \rho $}   & \footnotesize{$ u $}     & \footnotesize{$ p $}   & \footnotesize{$ Y_1 $} & \footnotesize{$\alpha_1$} & & & \footnotesize{$\gamma$  }\\\hline
   \footnotesize{Fluid 1}   &   \footnotesize{$ 10 $}   & \footnotesize{$ 50.0 $  } & \footnotesize{$ 1.1\cdot10^5  $} & \footnotesize{$ 1.0  $} & \footnotesize{$1.0$  } &&\footnotesize{Fluid 1}  &  \footnotesize{$1.4$  }  \\
   \footnotesize{Fluid 2}   &   \footnotesize{$ 1.0 $}   & \footnotesize{$ 50.0 $  } & \footnotesize{$ 1.0\cdot10^5  $} & \footnotesize{$ 0.0  $} & \footnotesize{$0.0$  } &&\footnotesize{Fluid 2}  &  \footnotesize{$1.1$  }  \\
   \bottomrule
\end{tabu}
\end{table}

\begin{figure}[h!]
\centering
\includegraphics[width=0.7\textwidth]{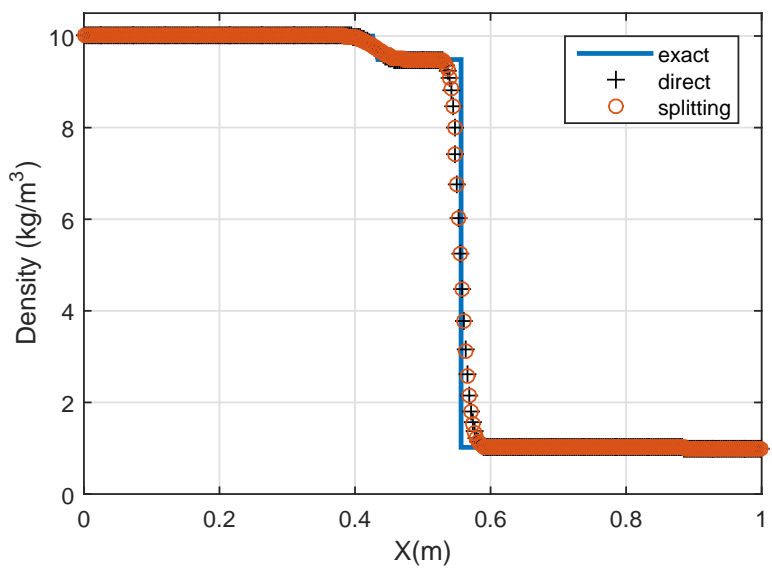}
\caption{Two-pressure jump problem - density profile - Exact solution ``{\color{blue}-}'', splitting approach ``{\color{Orange}$\boldsymbol{\circ}$}'' and direct approach ``$\bold{+}$"at $t=1.0$~ms.}
\label{Pressure_jump_nr8_density}
\end{figure}

\begin{figure}[h!]
\centering
\includegraphics[width=0.7\textwidth]{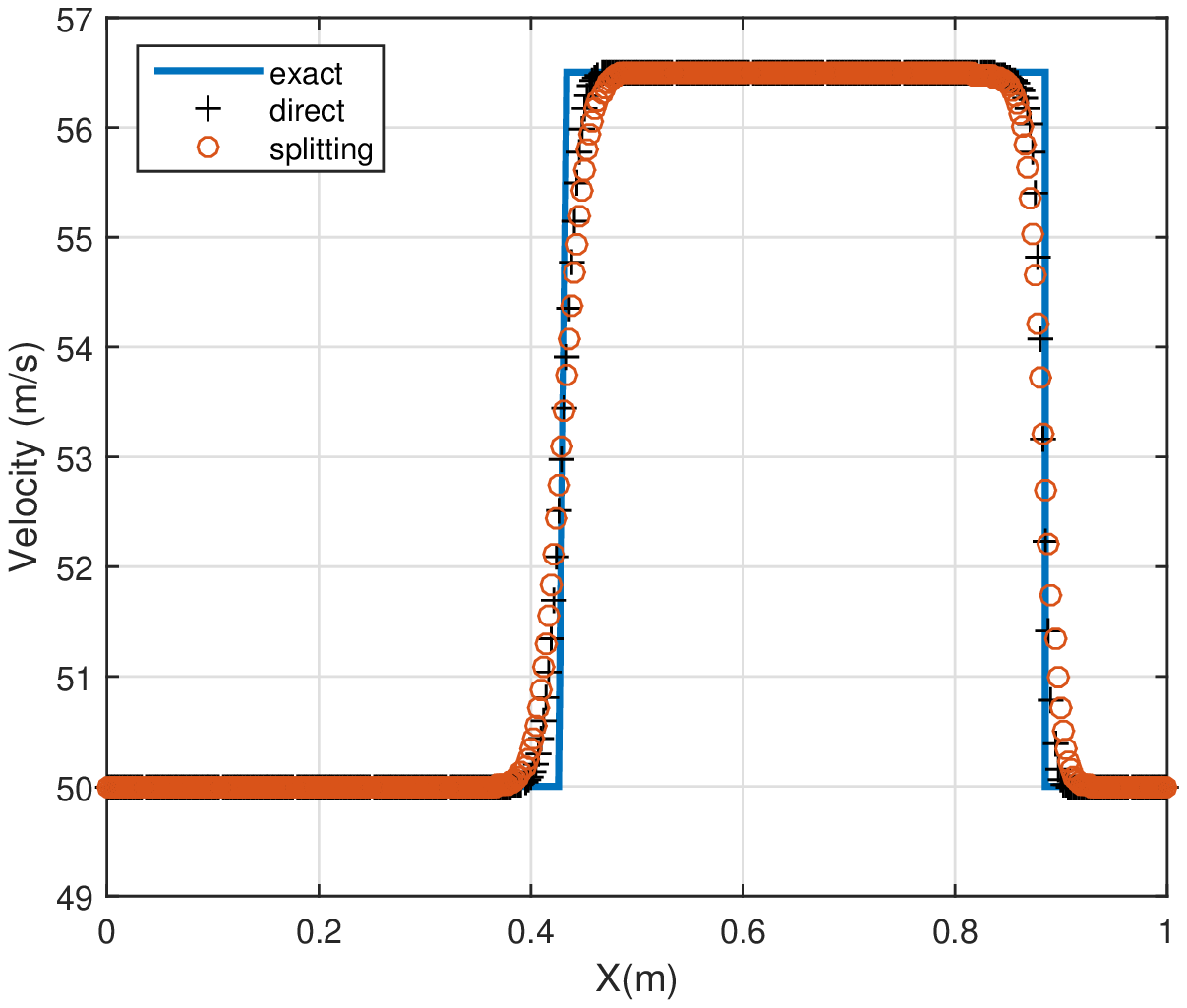}
\caption{Two-pressure jump problem - velocity profile- Exact solution ``{\color{blue}-}'', splitting approach ``{\color{Orange}$\boldsymbol{\circ}$}'' and direct approach ``$\bold{+}$"at $t=1.0$~ms.}
\label{Pressure_jump_nr8_velocity}
\end{figure}

\begin{figure}[h!]
\centering
\includegraphics[width=0.7\textwidth]{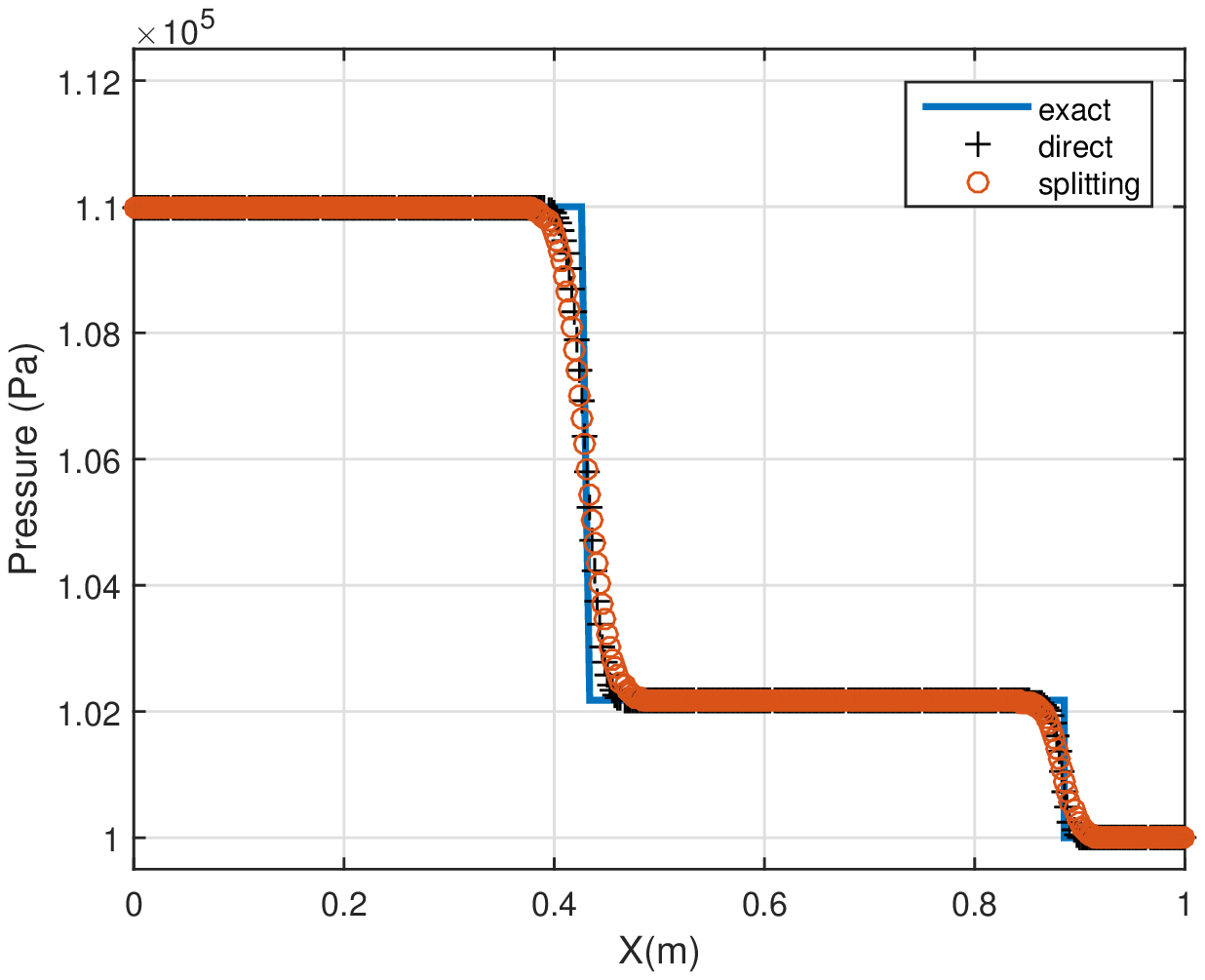}
\caption{Two-pressure jump problem - pressure profile - Exact solution ``{\color{blue}-}'', splitting approach ``{\color{Orange}$\boldsymbol{\circ}$}'' and direct approach ``$\bold{+}$"at $t=1.0$~ms.}
\label{Pressure_jump_nr8_pressure}
\end{figure}

\begin{figure}[h!]
\centering
\includegraphics[width=0.7\textwidth]{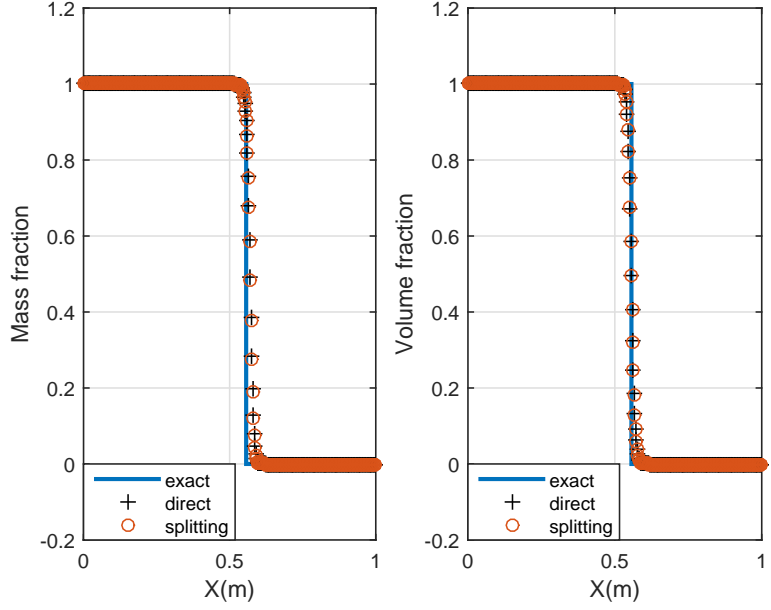}
\caption{Two-pressure jump problem - mass and volume fraction profile - Exact solution ``{\color{blue}-}'', splitting approach ``{\color{Orange}$\boldsymbol{\circ}$}'' and direct approach ``$\bold{+}$"at $t=1.0$~ms.}
\label{Pressure_jump_nr8_fractions}
\end{figure}

\begin{figure}[h!]
\centering
\includegraphics[width=1.0\textwidth]{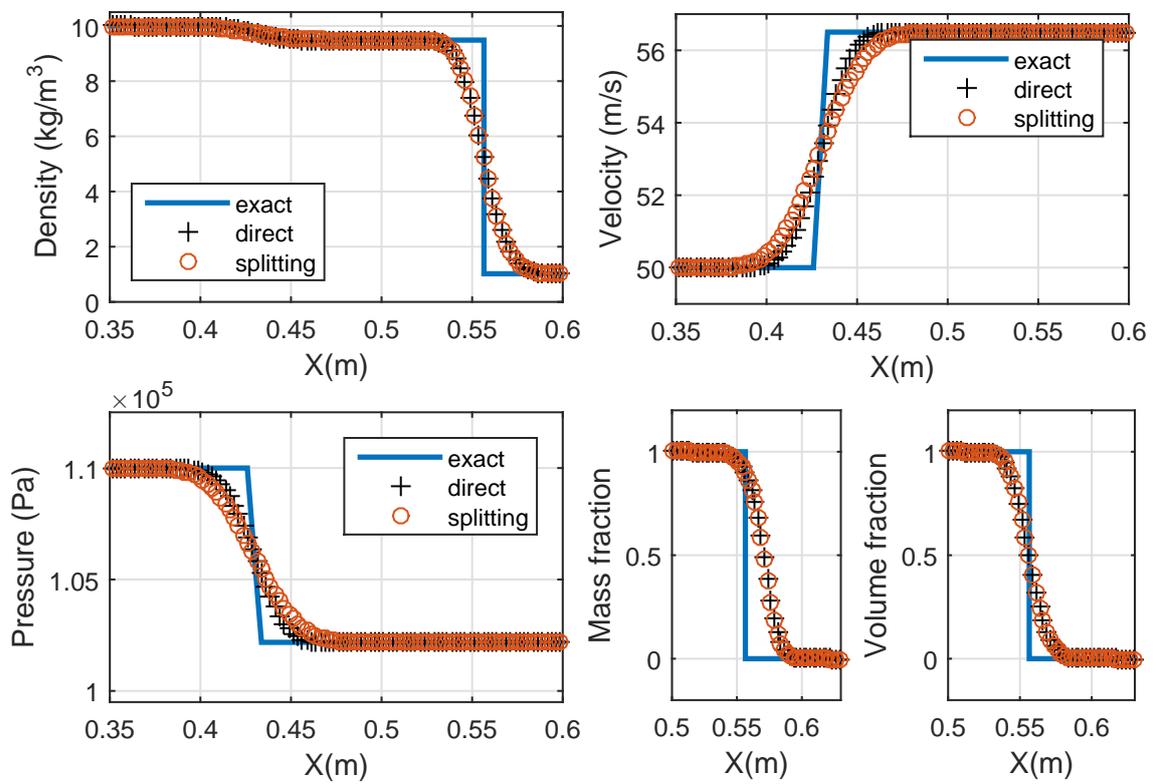}
\caption{Two-pressure jump problem - zoom - Exact solution ``{\color{blue}-}'', splitting approach ``{\color{Orange}$\boldsymbol{\circ}$}'' and direct approach ``$\bold{+}$"at $t=1.0$~ms.}
\label{Pressure_jump_nr8_zoom}
\end{figure}

The initial values and material properties are given in Table \ref{init pressure jump}. Also here the SG EOS reduces to the PG EOS. The results are obtained at time $t=1.0$~ms with $N=400$ cells for the Courant number of $\mathcal{C}=0.95$. The distributions of the primitive variables at $t=1.0$~ms are depicted in Figures \ref{Pressure_jump_nr8_density}-\ref{Pressure_jump_nr8_zoom} \added[id=Revs]{and the convergence rates are listed in Table \ref{convratesnr8}.}


\begin{table}[h!]
\caption{The $L_1$-convergence rates for the two-pressure jump problem. The convergence rates are computed as $c_N=\log(e_N/e_{2N})/\log(2)$. The errors are given by $e_N=\|s_N-s_{\text{exact}}\|_{L_1}$, where $s_N$ is the solution computed with $N$ grid points, $s_{\text{exact}}$ the exact solution, and $\|\cdot\|_{L_1}$ the standard $L_1$-norm.}\label{convratesnr8}
\begin{tabu} to \textwidth {X X X X X X X}
   \toprule
   \multicolumn{2}{c}{\footnotesize{Convergence rates}}& \multicolumn{5}{c}{\footnotesize{Physical quantity}}   \\
   \cline{1-7}
     \footnotesize{Splitting}                   & \footnotesize{approach}        & \footnotesize{$ \rho $}   & \footnotesize{$ u $}     & \footnotesize{$ p $}   & \footnotesize{$ Y_1 $} & \footnotesize{$\alpha_1$} \\ \hline
   \multicolumn{2}{c}{\footnotesize{$c_{40}$}}   &   \footnotesize{$ 0.43 $}   & \footnotesize{$ 0.69 $  } & \footnotesize{$ 0.65 $} & \footnotesize{$ 0.88  $} & \footnotesize{$0.41$  } \\
   \multicolumn{2}{c}{\footnotesize{$c_{80}$}}  &    \footnotesize{$ 0.54 $}   & \footnotesize{$ 0.59 $  } & \footnotesize{$ 0.49 $} & \footnotesize{$ 0.31  $} & \footnotesize{$0.56$  } \\
   \multicolumn{2}{c}{\footnotesize{$c_{160}$}}  &     \footnotesize{$ 0.50 $}   & \footnotesize{$ 0.65 $  } & \footnotesize{$ 0.59 $} & \footnotesize{$ 0.51  $} & \footnotesize{$0.50$  } \\
   \multicolumn{2}{c}{\footnotesize{$c_{320}$}}  &     \footnotesize{$ 0.50$}   & \footnotesize{$ 0.58 $  } & \footnotesize{$ 0.54 $} & \footnotesize{$ 0.50  $} & \footnotesize{$0.50$  } \\
   \multicolumn{2}{c}{\footnotesize{$c_{640}$}}   &     \footnotesize{$ 0.50 $}   & \footnotesize{$ 0.59 $  } & \footnotesize{$ 0.56 $} & \footnotesize{$ 0.50  $} & \footnotesize{$0.50$  } \\
   \cline{1-7}
     \footnotesize{Direct}                   & \footnotesize{approach}           & \footnotesize{$ \rho $}   & \footnotesize{$ u $}     & \footnotesize{$ p $}   & \footnotesize{$ Y_1 $} & \footnotesize{$\alpha_1$} \\ \hline
   \multicolumn{2}{c}{\footnotesize{$c_{40}$}}   &   \footnotesize{$ 0.42 $}   & \footnotesize{$ 0.70 $  } & \footnotesize{$ 0.69 $} & \footnotesize{$ 0.86  $} & \footnotesize{$0.40$  } \\
   \multicolumn{2}{c}{\footnotesize{$c_{80}$}}  &    \footnotesize{$ 0.54 $}   & \footnotesize{$ 0.56 $  } & \footnotesize{$ 0.46 $} & \footnotesize{$ 0.30  $} & \footnotesize{$0.56$  } \\
   \multicolumn{2}{c}{\footnotesize{$c_{160}$}}  &     \footnotesize{$ 0.50 $}   & \footnotesize{$ 0.69 $  } & \footnotesize{$ 0.61 $} & \footnotesize{$ 0.51  $} & \footnotesize{$0.49$  } \\
   \multicolumn{2}{c}{\footnotesize{$c_{320}$}}  &     \footnotesize{$ 0.50$}   & \footnotesize{$ 0.59 $  } & \footnotesize{$ 0.54 $} & \footnotesize{$ 0.50  $} & \footnotesize{$0.50$  } \\
   \multicolumn{2}{c}{\footnotesize{$c_{640}$}}   &     \footnotesize{$ 0.50 $}   & \footnotesize{$ 0.62 $  } & \footnotesize{$ 0.58 $} & \footnotesize{$ 0.50  $} & \footnotesize{$0.50$  } \\
      \bottomrule
\end{tabu}
\end{table}

Again, the results obtained with the proposed method are very similar to the ones obtained with the unsplit approach from \cite{daude2014numerical}. The location of the shock wave is accurately captured with both methods, also in the zoom (Figure \ref{Pressure_jump_nr8_zoom}) no significant difference is visible. Also for this test case, the newly proposed method takes larger time steps ($149$ time steps) than the direct approach from \cite{daude2014numerical} ($166$ time steps). \added[id=Rev.1]{The CPU time is $0.26$s and $0.69$s for the splitting approach and the direct approach respectively (averaged over $500$ runs on an i5 processor).} \added[id=Authors]{Again, both methods show similar convergence rates, see Table \ref{convratesnr8}.}
\clearpage

\subsection{No-reflection problem}


The third test we perform is the so-called no-reflection problem, which is also assessed in \cite{kreeft2010new}. In this test case, the right state is initially at rest and the left state moves towards the right state. The density and pressure of the left state are high compared with the right state. This will cause the two-fluid interface and the shock wave to move rightwards. The initial conditions \added[id=Authors]{of the perfect gases} are chosen such that no reflection wave occurs.
\begin{table}[h!]
\caption{Initial values and material properties for the no-reflection problem.}\label{init no refl}
\begin{tabu} to \textwidth {X X X X X X X X X X}
   \toprule
   \multicolumn{3}{c}{\footnotesize{(a) Initial values}} & &&&& \multicolumn{2}{c}{\footnotesize{(b) Material properties}}\\
   \cline{1-6} \cline{8-9}
   & \footnotesize{$ \rho $}   & \footnotesize{$ u $}     & \footnotesize{$ p $}   & \footnotesize{$ Y_1 $} & \footnotesize{$\alpha_1$} & & & \footnotesize{$\gamma$  }\\\hline
   \footnotesize{Fluid 1}   &   \footnotesize{$3.1748 $}   & \footnotesize{$  9.4350 $  } & \footnotesize{$ 100  $} & \footnotesize{$ 1.0  $} & \footnotesize{$1.0$  } &&\footnotesize{Fluid 1}  &  \footnotesize{$1.667$  }  \\
   \footnotesize{Fluid 2}   &   \footnotesize{$ 1.0 $}   & \footnotesize{$ 0.0 $  } & \footnotesize{$ 1.0  $} & \footnotesize{$ 0.0  $} & \footnotesize{$0.0$  } &&\footnotesize{Fluid 2}  &  \footnotesize{$1.2$  } \\
   \bottomrule
\end{tabu}
\end{table}

\begin{figure}[h!]
\centering
\includegraphics[width=0.6\textwidth]{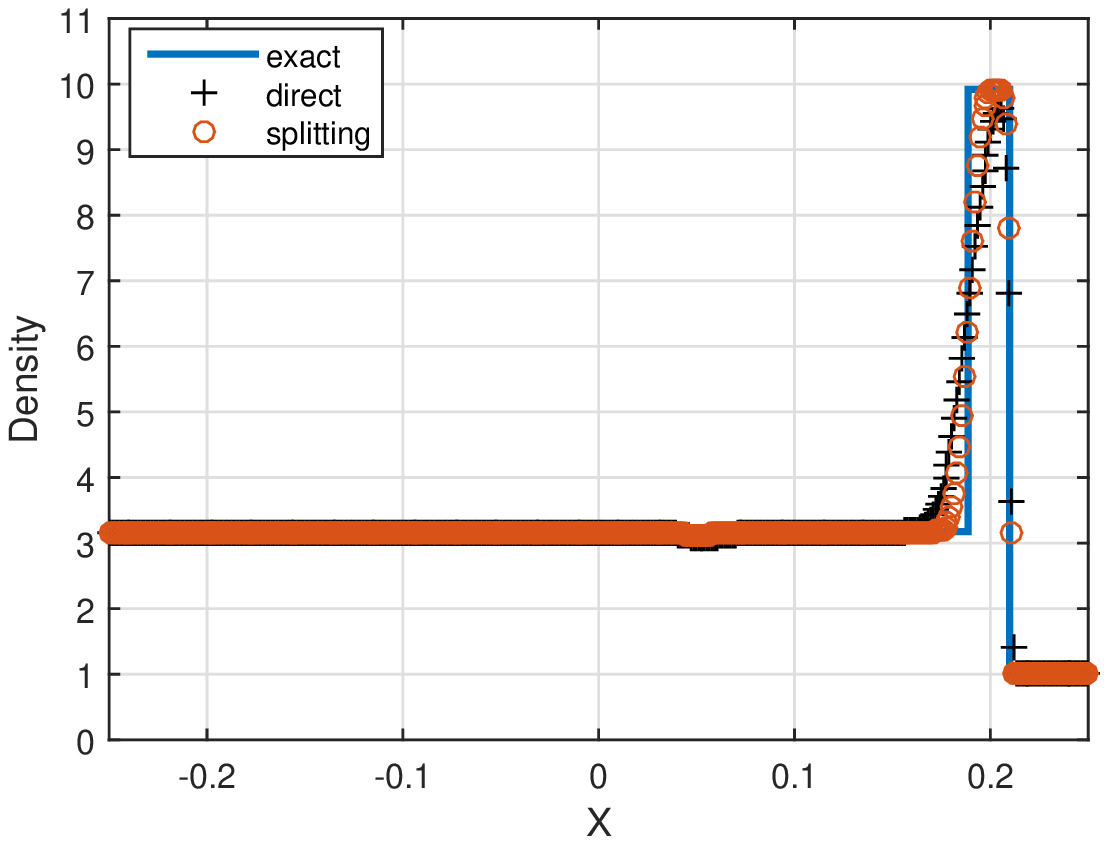}
\caption{No-reflection problem - density profile - Exact solution ``{\color{blue}-}'', splitting approach ``{\color{Orange}$\boldsymbol{\circ}$}'' and direct approach ``$\bold{+}$"at $t=0.02$.}
\label{norefl_density}
\end{figure}

\begin{figure}[h!]
\centering
\includegraphics[width=0.6\textwidth]{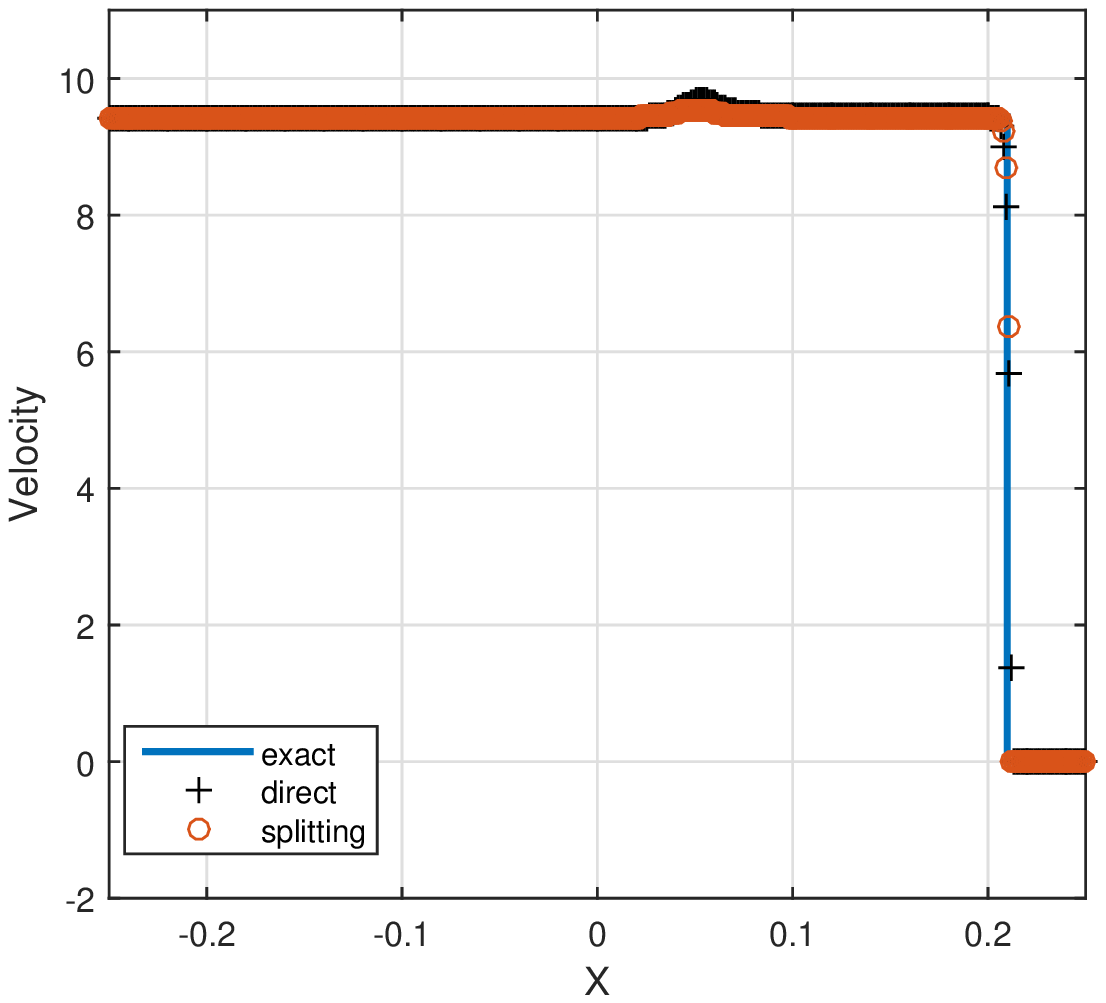}
\caption{No-reflection problem - velocity profile - Exact solution ``{\color{blue}-}'', splitting approach ``{\color{Orange}$\boldsymbol{\circ}$}'' and direct approach ``$\bold{+}$"at $t=0.02$.}
\label{norefl_velocity}
\end{figure}

\begin{figure}[h!]
\centering
\includegraphics[width=0.6\textwidth]{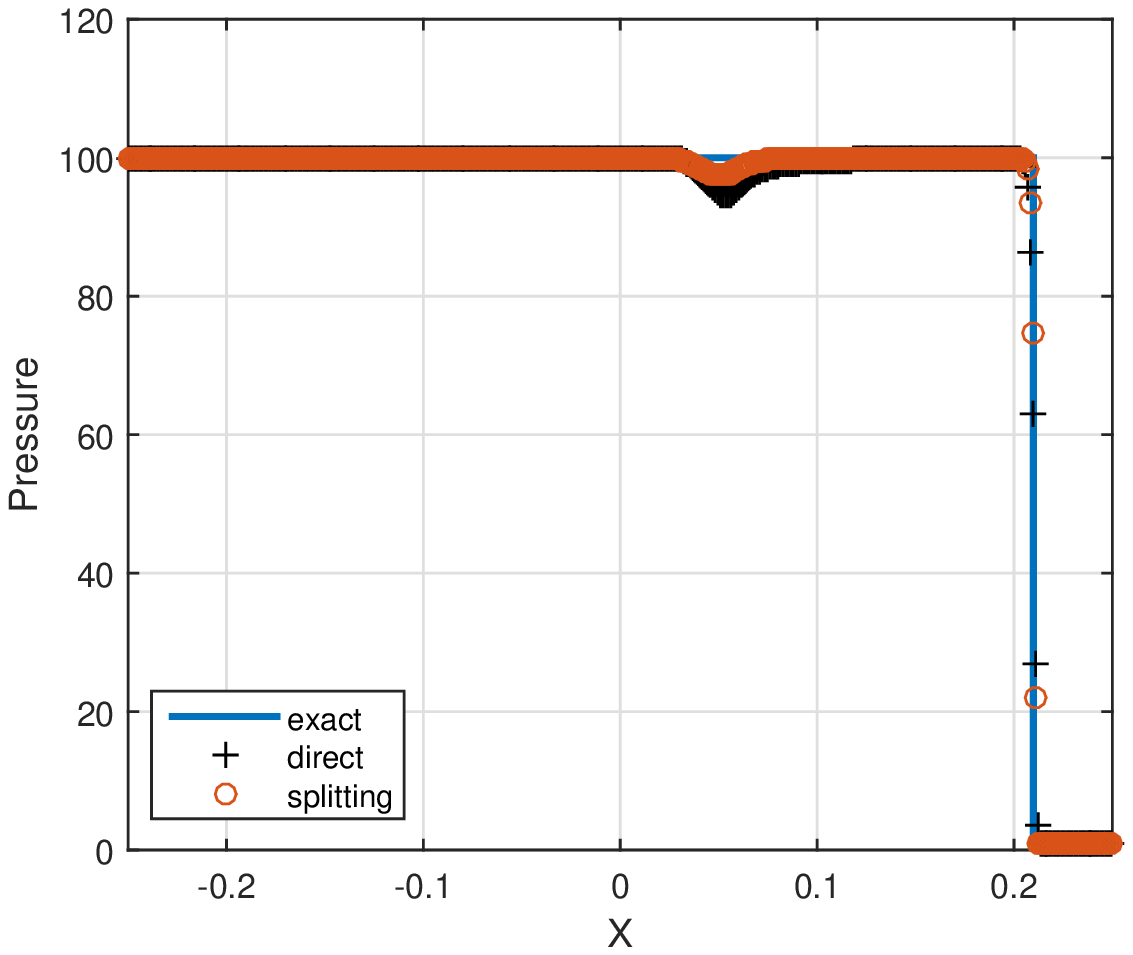}
\caption{No-reflection problem - pressure profile - Exact solution ``{\color{blue}-}'', splitting approach ``{\color{Orange}$\boldsymbol{\circ}$}'' and direct approach ``$\bold{+}$"at $t=0.02$.}
\label{norefl_pressure}
\end{figure}

\begin{figure}[h!]
\centering
\includegraphics[width=0.6\textwidth]{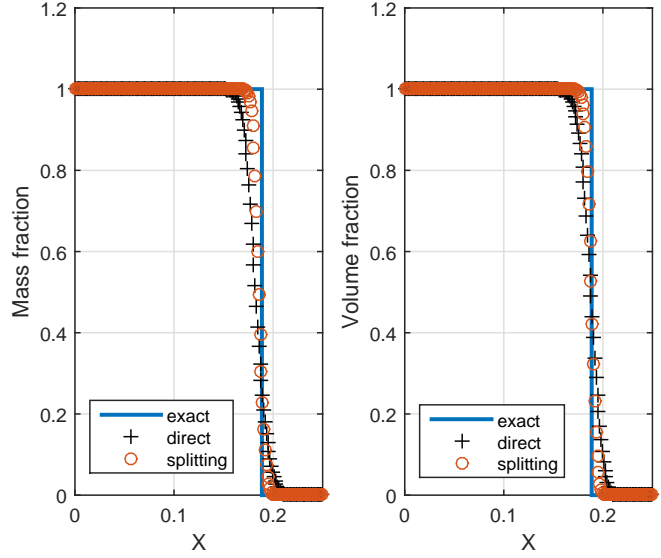}
\caption{No-reflection problem - mass and volume fraction profile - Exact solution ``{\color{blue}-}'', splitting approach ``{\color{Orange}$\boldsymbol{\circ}$}'' and direct approach ``$\bold{+}$"at $t=0.02$.}
\label{norefl_fractions}
\end{figure}

\begin{figure}[h!]
\centering
\includegraphics[width=1.0\textwidth]{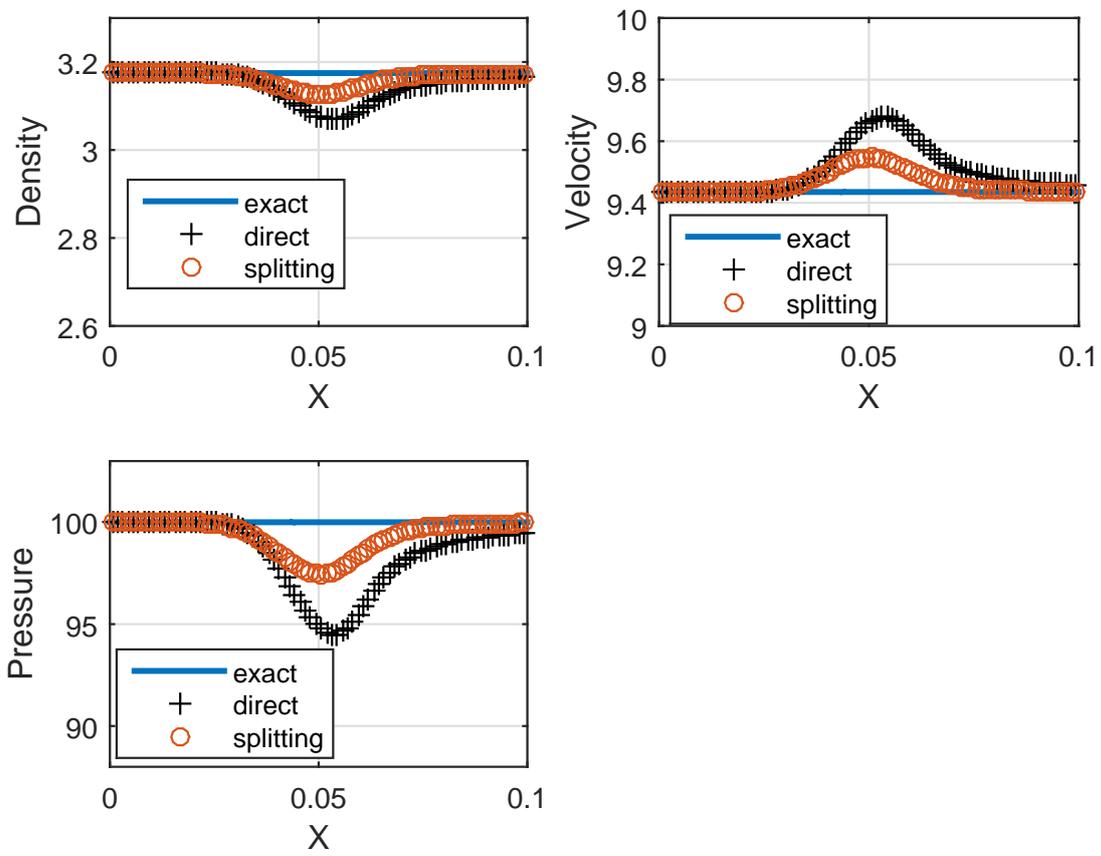}
\caption{No-reflection problem - zoom at bumps - Exact solution ``{\color{blue}-}'', splitting approach ``{\color{Orange}$\boldsymbol{\circ}$}'' and direct approach ``$\bold{+}$"at $t=0.02$.}
\label{norefl_zoom_bump}
\end{figure}

\begin{figure}[h!]
\centering
\includegraphics[width=1.0\textwidth]{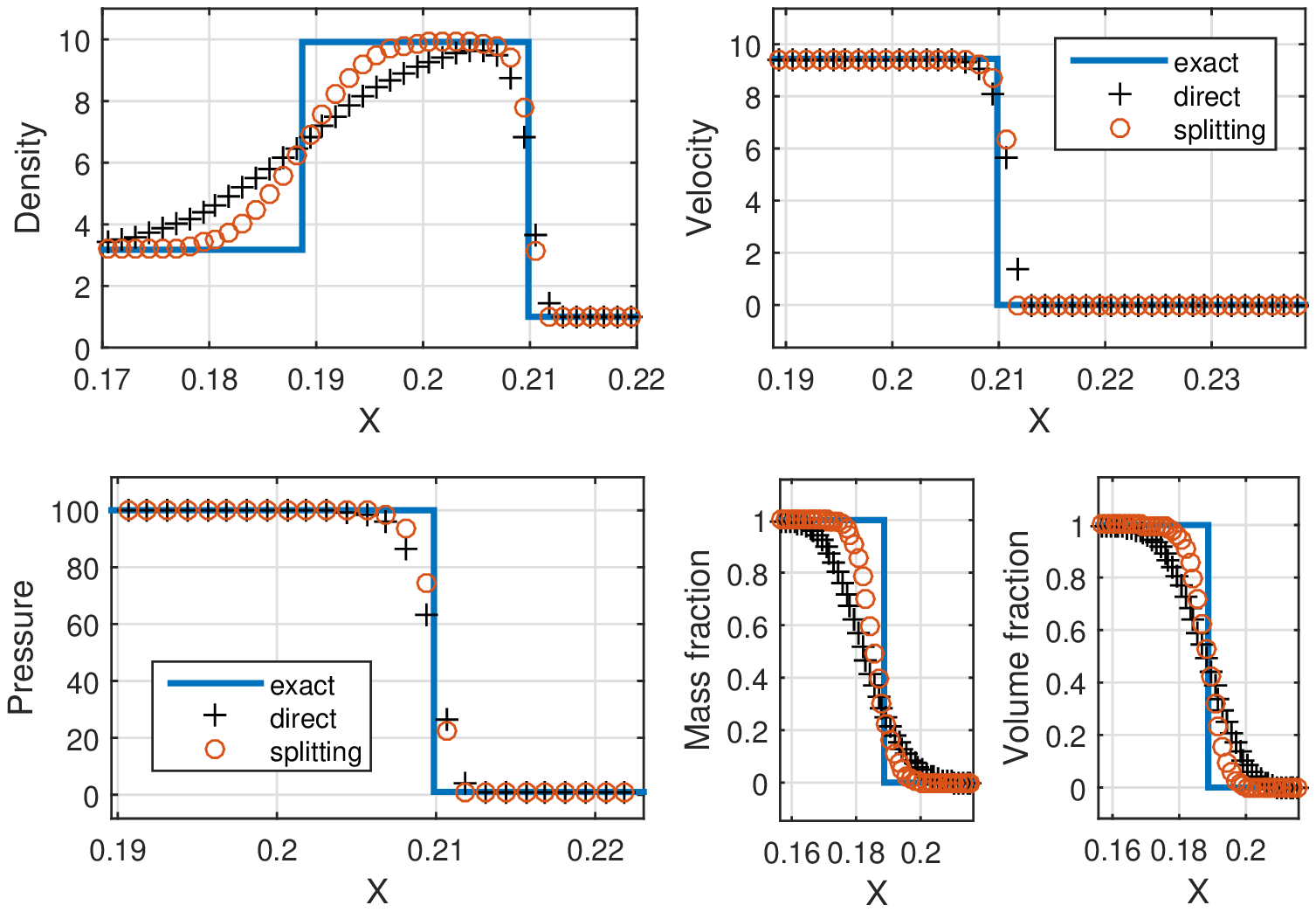}
\caption{No-reflection problem - zoom at shock wave - Exact solution ``{\color{blue}-}'', splitting approach ``{\color{Orange}$\boldsymbol{\circ}$}'' and direct approach ``$\bold{+}$"at $t=0.02$.}
\label{norefl_zoom_shock}
\end{figure}

The initial values and material properties are given in Table \ref{init no refl}. The results are obtained at time $t=0.02$ with $N=400$ cells \added[id=Authors]{with a CFL number of $\mathcal{C}=0.95$.} The results are visualized in the Figures \ref{norefl_density}-\ref{norefl_zoom_shock} \added[id=Revs]{and the convergence rates are listed in Table \ref{convratesnr5}.}
\\
\begin{table}[h!]
\caption{The $L_1$-convergence rates for the no-reflection problem. The errors are computed as $e_N=\|s_N-s_{\text{exact}}\|_{L_1}$, where $s_N$ is the solution computed with $N$ grid points, $s_{\text{exact}}$ the exact solution, and $\|\cdot\|_{L_1}$ the standard $L_1$-norm.}\label{convratesnr5}
\begin{tabu} to \textwidth {X X X X X X X}
   \toprule
   \multicolumn{2}{c}{\footnotesize{Fraction errors}}& \multicolumn{5}{c}{\footnotesize{Physical quantity}}   \\
   \cline{1-7}
     \footnotesize{Splitting}                   & \footnotesize{approach}        & \footnotesize{$ \rho $}   & \footnotesize{$ u $}     & \footnotesize{$ p $}   & \footnotesize{$ Y_1 $} & \footnotesize{$\alpha_1$} \\ \hline
   \multicolumn{2}{c}{\footnotesize{$c_{40}$}}   &   \footnotesize{$ 0.69 $}   & \footnotesize{$ 1.05 $  } & \footnotesize{$ 1.04 $} & \footnotesize{$ 0.50  $} & \footnotesize{$0.45$  } \\
   \multicolumn{2}{c}{\footnotesize{$c_{80}$}}  &    \footnotesize{$ 0.82 $}   & \footnotesize{$ 1.33 $  } & \footnotesize{$ 1.22 $} & \footnotesize{$ 0.52  $} & \footnotesize{$0.46$  } \\
   \multicolumn{2}{c}{\footnotesize{$c_{160}$}}  &     \footnotesize{$ 0.57 $}   & \footnotesize{$ 0.83 $  } & \footnotesize{$ 0.90 $} & \footnotesize{$ 0.52  $} & \footnotesize{$0.45$  } \\
   \multicolumn{2}{c}{\footnotesize{$c_{320}$}}  &     \footnotesize{$ 0.57$}   & \footnotesize{$ 0.80 $  } & \footnotesize{$ 0.82 $} & \footnotesize{$ 0.43  $} & \footnotesize{$0.50$  } \\
   \multicolumn{2}{c}{\footnotesize{$c_{640}$}}   &     \footnotesize{$ 0.68 $}   & \footnotesize{$ 1.27 $  } & \footnotesize{$ 1.22 $} & \footnotesize{$ 0.43  $} & \footnotesize{$0.50$  } \\
   \cline{1-7}
     \footnotesize{Direct}                   & \footnotesize{approach}           & \footnotesize{$ \rho $}   & \footnotesize{$ u $}     & \footnotesize{$ p $}   & \footnotesize{$ Y_1 $} & \footnotesize{$\alpha_1$} \\ \hline
   \multicolumn{2}{c}{\footnotesize{$c_{40}$}}   &   \footnotesize{$ 0.44 $}   & \footnotesize{$ 1.03 $  } & \footnotesize{$ 0.86 $} & \footnotesize{$ 0.42  $} & \footnotesize{$0.42$  } \\
   \multicolumn{2}{c}{\footnotesize{$c_{80}$}}  &    \footnotesize{$ 0.54 $}   & \footnotesize{$ 1.17 $  } & \footnotesize{$ 0.93 $} & \footnotesize{$ 0.42  $} & \footnotesize{$0.38$  } \\
   \multicolumn{2}{c}{\footnotesize{$c_{160}$}}  &     \footnotesize{$ 0.50 $}   & \footnotesize{$ 0.87 $  } & \footnotesize{$ 0.90 $} & \footnotesize{$ 0.45  $} & \footnotesize{$0.39$  } \\
   \multicolumn{2}{c}{\footnotesize{$c_{320}$}}  &     \footnotesize{$ 0.51$}   & \footnotesize{$ 0.81 $  } & \footnotesize{$ 0.91 $} & \footnotesize{$ 0.41  $} & \footnotesize{$0.42$  } \\
   \multicolumn{2}{c}{\footnotesize{$c_{640}$}}   &     \footnotesize{$ 0.56 $}   & \footnotesize{$ 1.20 $  } & \footnotesize{$ 1.06 $} & \footnotesize{$ 0.44  $} & \footnotesize{$0.45$  } \\
      \bottomrule
\end{tabu}
\end{table}
The \added[id=Authors]{location} of the contact discontinuity \added[id=Authors]{is} satisfactorily retrieved with \added[id=Authors]{both} methods. \added[id=Authors]{A} small reflected wave is visible at around $x=0.05$, which is \added[id=Authors]{weaker} for the splitting-based scheme (see Figure \ref{norefl_zoom_bump}). For both methods it vanishes when refining the grid. The \added[id=Authors]{shock wave is well retrieved with both methods. The newly proposed method seems to be less diffusive than the direct approach (see Figure \ref{norefl_zoom_shock})}. Again, the newly proposed method takes larger time steps \added[id=Authors]{($169$ time steps)} than the direct approach from \cite{daude2014numerical} \added[id=Authors]{($285$ time steps)}. \added[id=Rev.1]{The CPU time is $0.25$s and $0.33$s for the splitting approach and the direct approach, respectively (averaged over $500$ runs on an i5 processor).}
\clearpage

\subsection{Water-air mixture problem}

\begin{figure}[h!]
\centering
\includegraphics[width=0.7\textwidth]{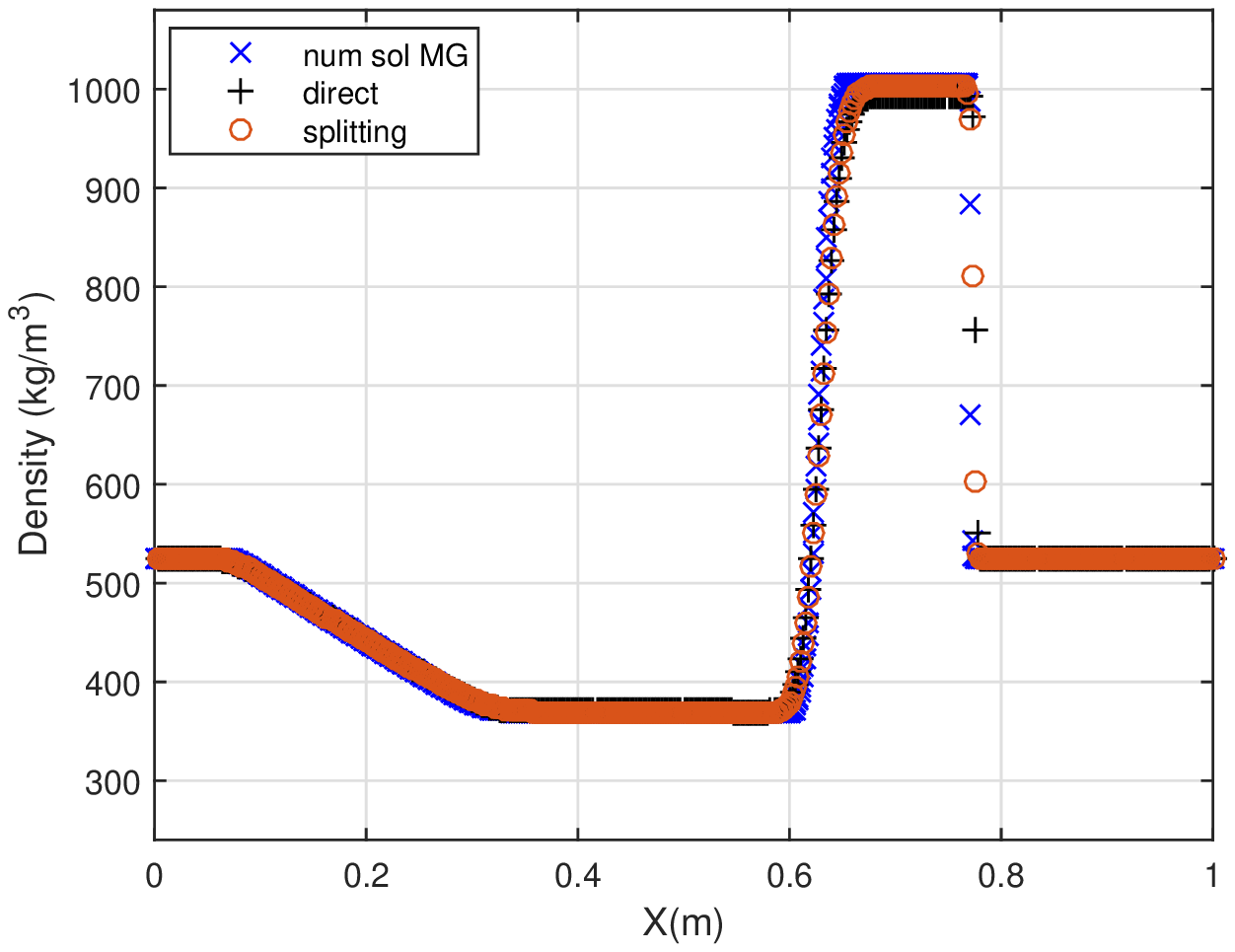}
            \caption{Water-air mixture problem - Density profile - Numerical solution from \cite{murrone2005five} ``{\color{blue}x}'', splitting approach ``{\color{Orange}$\boldsymbol{\circ}$}'' and direct approach ``$\bold{+}$"at $t=200 \mu$ s.}
            \label{mixture_all nr10_density}
\end{figure}
\begin{figure}[h!]
\centering
\includegraphics[width=0.7\textwidth]{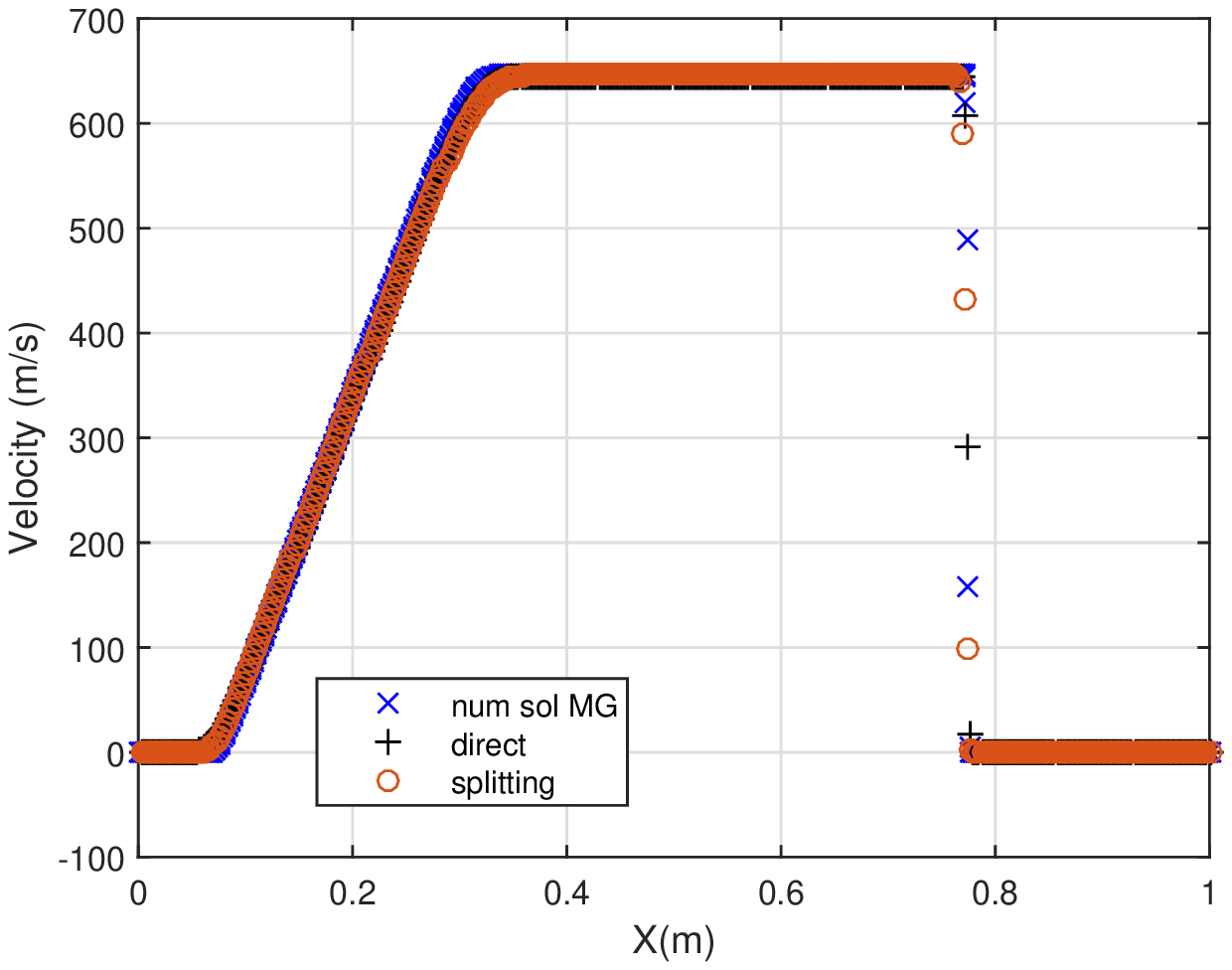}
            \caption{Water-air mixture problem - Velocity profile - Numerical solution from \cite{murrone2005five} ``{\color{blue}x}'', splitting approach ``{\color{Orange}$\boldsymbol{\circ}$}'' and direct approach ``$\bold{+}$"at $t=200 \mu$ s.}
            \label{mixture_all nr10_velocity}
\end{figure}
\begin{figure}[h!]
\centering
\includegraphics[width=0.7\textwidth]{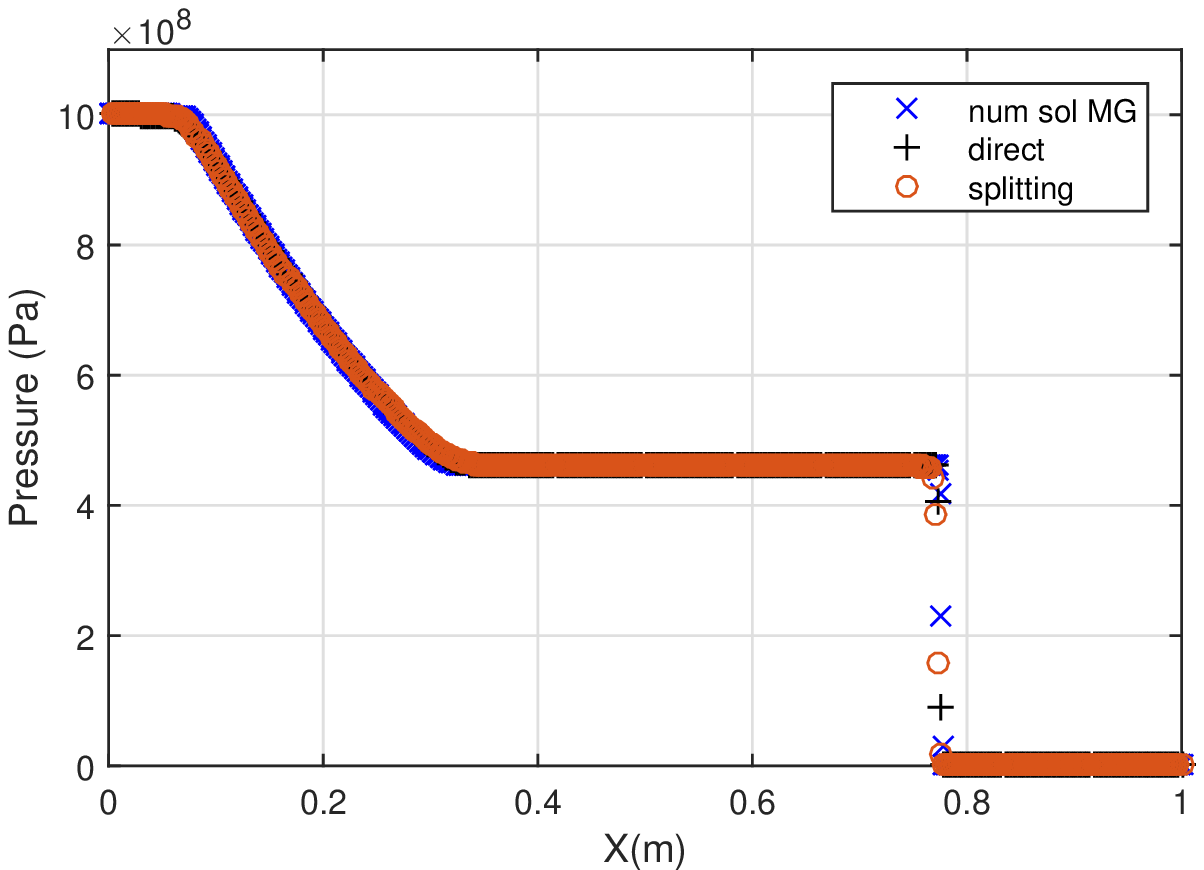}
            \caption{Water-air mixture problem - Pressure profile - Numerical solution from \cite{murrone2005five} ``{\color{blue}x}'', splitting approach ``{\color{Orange}$\boldsymbol{\circ}$}'' and direct approach ``$\bold{+}$"at $t=200 \mu$ s.}
            \label{mixture_all nr10_pressure}
\end{figure}
\begin{figure}[h!]
\centering
\includegraphics[width=0.7\textwidth]{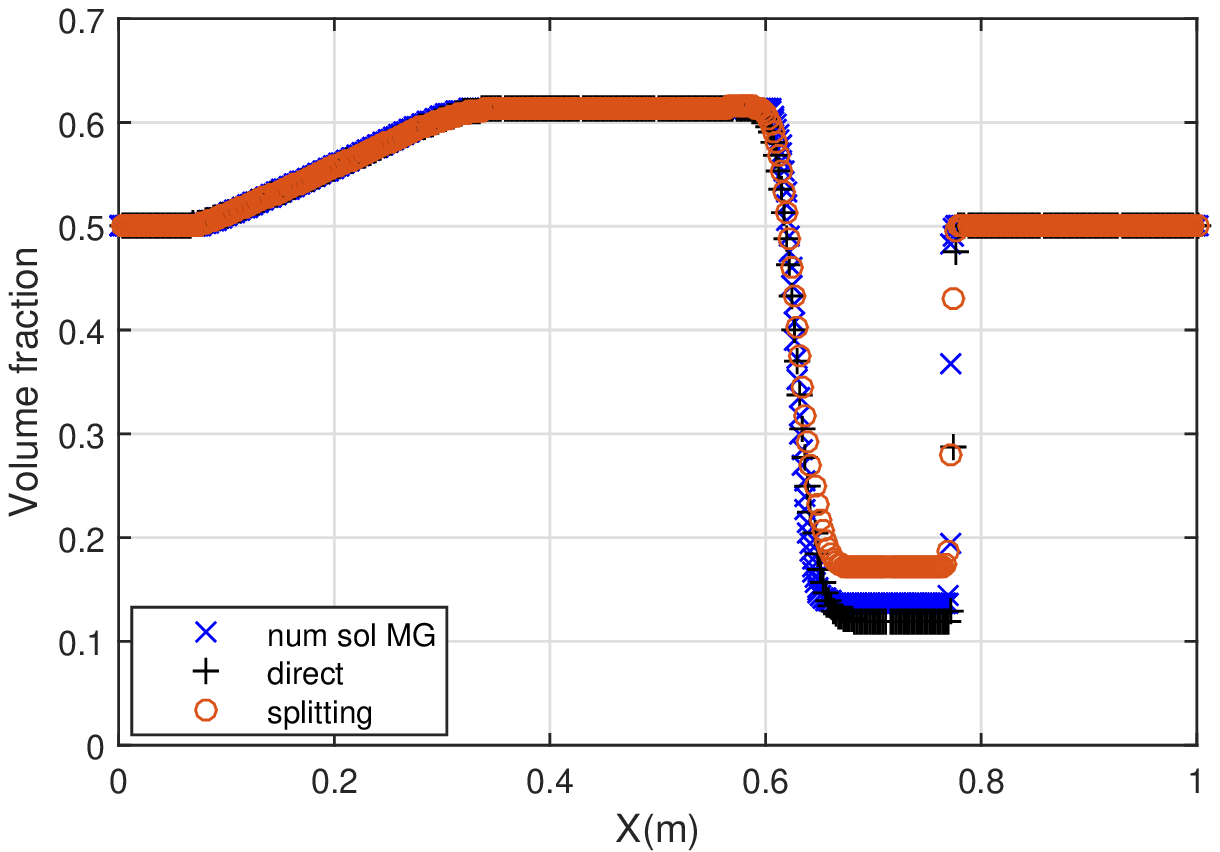}
            \caption{Water-air mixture problem - Volume fraction profile - Numerical solution from \cite{murrone2005five} ``{\color{blue}x}'', splitting approach ``{\color{Orange}$\boldsymbol{\circ}$}'' and direct approach ``$\bold{+}$"at $t=200 \mu$ s.}
            \label{mixture_all nr10_alpha}
\end{figure}
\begin{figure}[h!]
\centering
\includegraphics[width=1.0\textwidth]{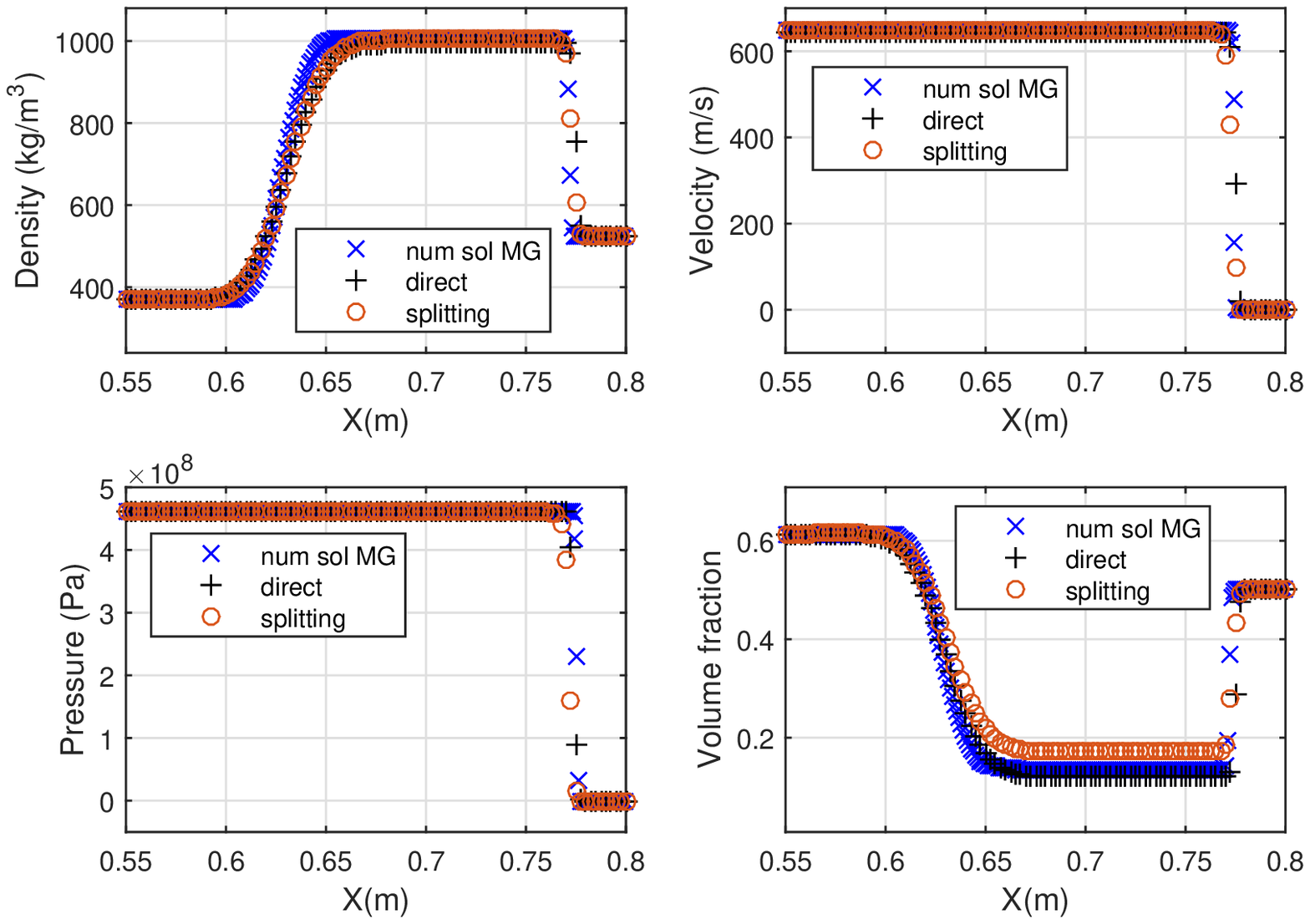}
            \caption{Water-air mixture problem - zoom - Numerical solution from \cite{murrone2005five} ``{\color{blue}x}'', splitting approach ``{\color{Orange}$\boldsymbol{\circ}$}'' and direct approach ``$\bold{+}$"at $t=200 \mu$ s.}
            \label{mixture_all nr10_zoom}
\end{figure}

In this shock tube test we consider a water-air mixture problem. This test case has been considered by Murrone and Guillard \cite{murrone2005five} and by Kreeft and Koren \cite{kreeft2010new}. In contrast to the previous test cases, the shock tube is now filled with a mixture of water and air ($0< Y_1, \alpha_1 < 1$) \added[id=Authors]{and stiffened gases are considered.} Both mixture states are initially at rest and the initial pressure ratio is $10^4$.\\
\begin{table}[h!]
\caption{Initial values and material properties for the water-air mixture problem. The dimensions of the quantities $\rho, u$ and $p$ are kg m$^{-3}$, m s$^{-1}$ and Pa respectively.}\label{init water air}
\begin{tabu} to \textwidth {X X X X X X}
   \toprule
   & \footnotesize{$ \rho $}   & \footnotesize{$ u $}     & \footnotesize{$ p $}   & \footnotesize{$ Y_1 $} & \footnotesize{$\alpha_1$}\\\hline
   \footnotesize{Left chamber}   &   \footnotesize{$ 525 $}   & \footnotesize{$ 0.0 $  } & \footnotesize{$ 10^9 $} & \footnotesize{$ 0.0476  $} & \footnotesize{$0.5$  }\\
   \footnotesize{Right chamber}   &   \footnotesize{$ 525 $}   & \footnotesize{$ 0.0 $  } & \footnotesize{$ 10^5 $} & \footnotesize{$ 0.9524  $} & \footnotesize{$0.5$  }\\
   \bottomrule
\end{tabu}
\end{table}
\begin{table}[h!]
\caption{Material properties for the water-air mixture problem. The dimensions of the quantities $\pi $ and $\eta$ are Pa and J kg$^{-1}$ respectively.}\label{init water air material}
\begin{tabu} to \textwidth {X X X X}
   \toprule
                                  & \footnotesize{$\gamma$  }    & \footnotesize{$\pi$  }           &  \footnotesize{$\eta$  }  \\\hline
   \footnotesize{Fluid 1}    & \footnotesize{$1.4$  }       &  \footnotesize{$0.0$  }          &  \footnotesize{$0.0$  } \\
   \footnotesize{Fluid 2}   & \footnotesize{$4.4$  }       &  \footnotesize{$6\cdot10^8$  }   &  \footnotesize{$0.0$  } \\
   \bottomrule
\end{tabu}
\end{table}\\
The initial values and material properties are given in Table\added[id=Authors]{s} \ref{init water air} \added[id=Authors]{and} \ref{init water air material}. Numerical results are obtained at time $t=200 ~\mu$s with $N=400$ cells with CFL number $\mathcal{C}=0.95$. The results are visualized in the Figures \ref{mixture_all nr10_density}-\ref{mixture_all nr10_zoom}.
\\\\
\indent The numerical results are \added[id=Authors]{in good agreement} with numerical solutions from Murrone and Guillard \cite{murrone2005five}. The volume fraction distribution on the right side of the middle wave shows slightly different values for all three schemes. The numerical solution from \cite{murrone2005five} shows a slightly lower value compared with the splitting-based method and \added[id=Authors]{a} slightly higher \added[id=Authors]{value} than the HLLC-type scheme. This test case indicates that the proposed method can also deal with mixture problems. Also for this test case, the newly proposed method takes larger time steps ($179$ time steps) than the direct approach from \cite{daude2014numerical} (giving $193$ time steps). \added[id=Rev.1]{The CPU time is $0.25$s and $0.38$s for the splitting approach and the direct approach, respectively (averaged over $500$ runs on an i5 processor).}
\subsection{\added[id=Rev.3]{Two-phase cavitation problem}}


In this test case proposed by Saurel et al. \cite{saurel2008modelling} the tube is filled with water and its vapor at atmospheric pressure. Thus a mixture of the fluids is considered: initially the water (with density $\rho_2=1150 $ kg m$^{-3}$) contains a small portion of vapor $\alpha_1=10^{-2}$ (with density $\rho_1=0.63$ kg m$^{-3}$). An initial velocity discontinuity separates both states.\\
\begin{table}[h!]
\caption{Initial values for the two-phase cavitation problem. The dimensions of the quantities $\rho, u$ and $p$ respectively.}\label{init cavitation}
\begin{tabu} to \textwidth {X X X X X X}
   \toprule
   & \footnotesize{$ \rho $}   & \footnotesize{$ u $}     & \footnotesize{$ p $}   & \footnotesize{$ Y_1 $} & \footnotesize{$\alpha_1$}\\\hline
   \footnotesize{Left chamber}   &   \footnotesize{$ 1138.5063 $}   & \footnotesize{$ -2.0 $  } & \footnotesize{$ 10^5 $} & \footnotesize{$ 5.53356 \cdot 10^{-6}  $} & \footnotesize{$0.01$  } \\
   \footnotesize{Right chamber}   &   \footnotesize{$ 1138.5063 $}   & \footnotesize{$ 2.0 $  } & \footnotesize{$ 10^5 $} & \footnotesize{$ 5.53356 \cdot 10^{-6}  $} & \footnotesize{$0.01$  } \\
   \bottomrule
\end{tabu}
\end{table}
\begin{table}[h!]
\caption{Material properties for the two-phase cavitation problem. The dimensions of the quantities $\pi $ and $\eta$ are Pa and J kg$^{-1}$ respectively.}\label{init cavitation material}
\begin{tabu} to \textwidth {X X X X}
   \toprule
                                  & \footnotesize{$\gamma$  }    & \footnotesize{$\pi$  }   &  \footnotesize{$\eta$  }  \\\hline
   \footnotesize{Fluid 1}   & \footnotesize{$2.35$  }      &  \footnotesize{$10^9$  } &  \footnotesize{$-1167\cdot10^3$  } \\
   \footnotesize{Fluid 2}   & \footnotesize{$1.43$  }      &  \footnotesize{$0$  }    &  \footnotesize{$2030\cdot10^3$  } \\
   \bottomrule
\end{tabu}
\end{table}

\begin{figure}[h!]
\centering
\includegraphics[width=0.7\textwidth]{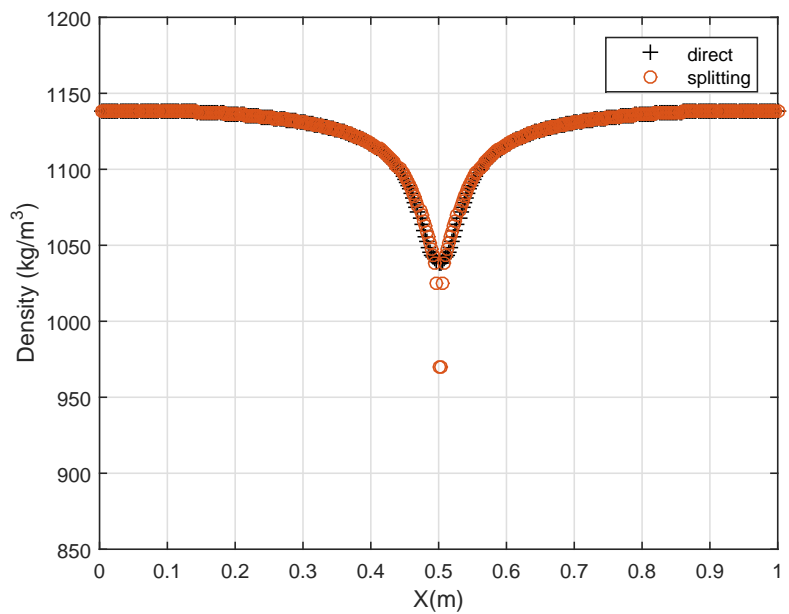}
            \caption{Two-phase cavitation problem problem - Density profile - Splitting approach ``{\color{Orange}$\boldsymbol{\circ}$}'' and direct approach ``$\bold{+}$"at $t=3.2$~ms.}
            \label{Two-phase cavitation problem nr-2_density}
\end{figure}

\begin{figure}[h!]
\centering
\includegraphics[width=0.7\textwidth]{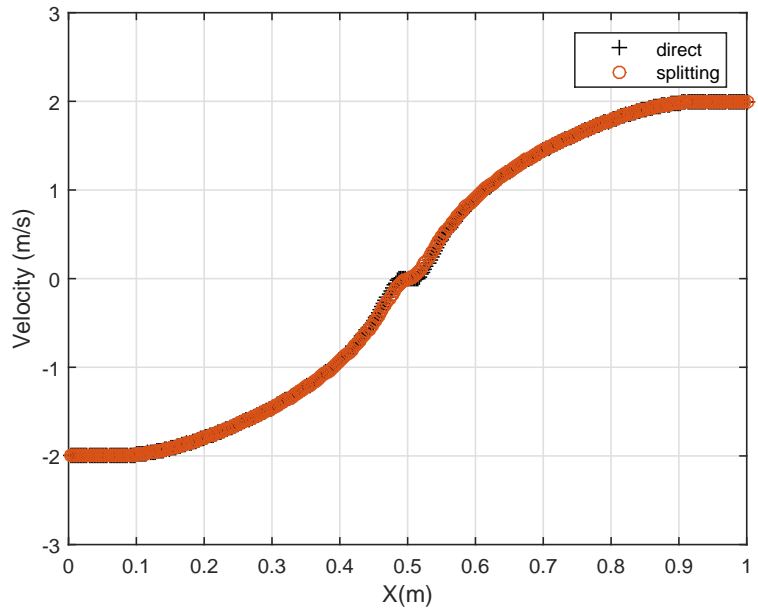}
            \caption{Two-phase cavitation problem problem - Velocity profile - Splitting approach ``{\color{Orange}$\boldsymbol{\circ}$}'' and direct approach ``$\bold{+}$"at $t=3.2$~ms.}
            \label{Two-phase cavitation problem nr-2_velocity}
\end{figure}

\begin{figure}[h!]
\centering
\includegraphics[width=0.7\textwidth]{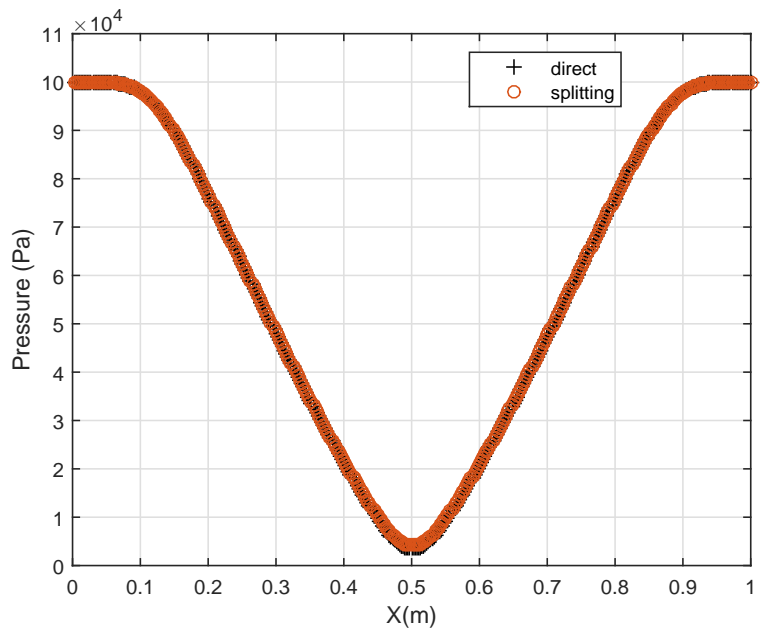}
            \caption{Two-phase cavitation problem problem - Pressure profile - Splitting approach ``{\color{Orange}$\boldsymbol{\circ}$}'' and direct approach ``$\bold{+}$"at $t=3.2$~ms.}
            \label{Two-phase cavitation problem nr-2_pressure}
\end{figure}

\begin{figure}[h!]
\centering
\includegraphics[width=0.7\textwidth]{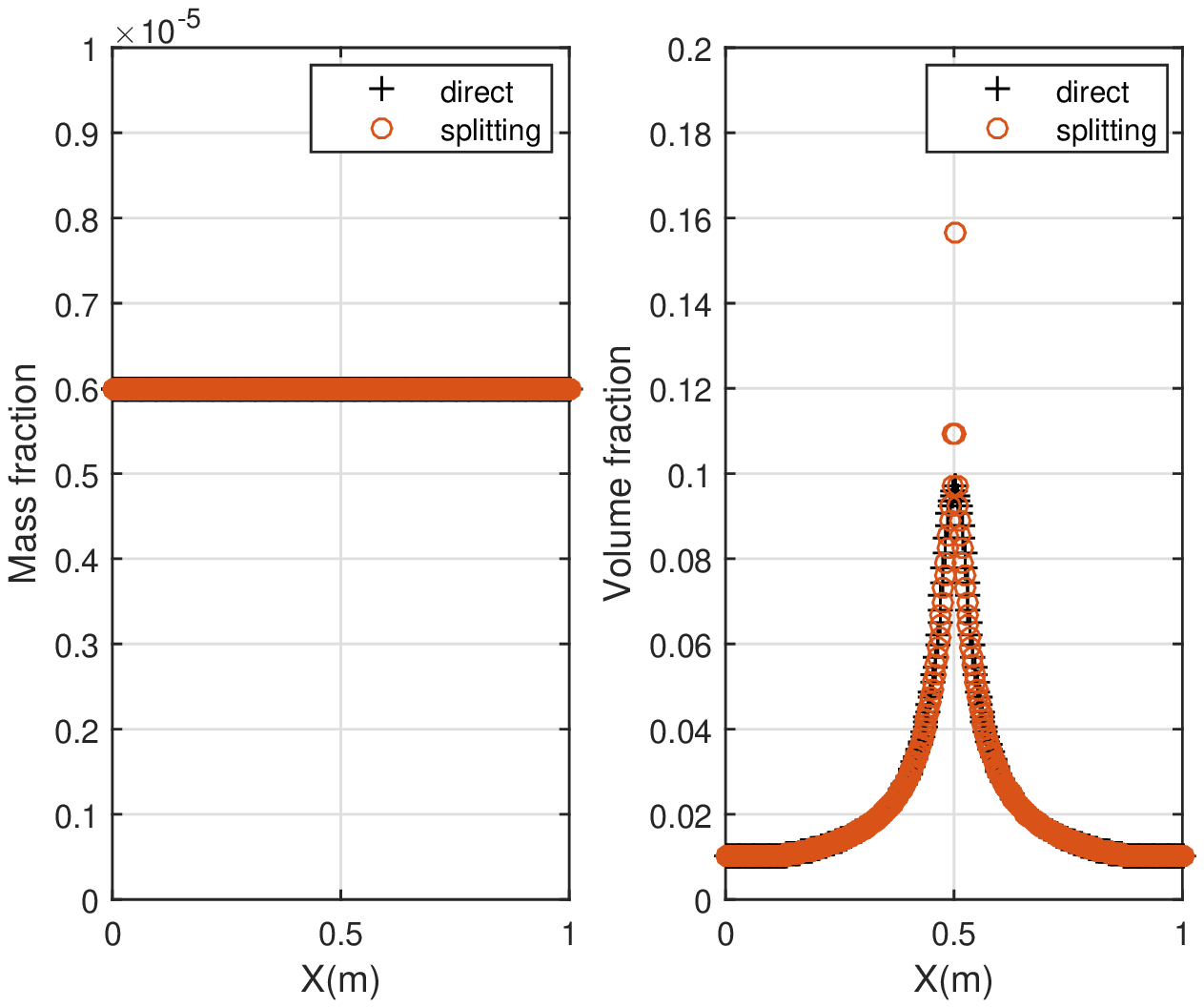}
            \caption{Two-phase cavitation problem problem - Mass and volume fraction profile - Splitting approach ``{\color{Orange}$\boldsymbol{\circ}$}'' and direct approach ``$\bold{+}$"at $t=3.2$~ms.}
            \label{Two-phase cavitation problem nr-2_fractions}
\end{figure}

The initial values and material properties are given in Tables \ref{init cavitation} and \ref{init cavitation material}. Numerical results are presented at time $t=3.2 ~m$s with $N=400$ cells. A smaller time step (CFL = $\mathcal{C}=0.01$) is used due to the strong rarefaction wave. The results are visualized in the Figures \ref{Two-phase cavitation problem nr-2_density}-\ref{Two-phase cavitation problem nr-2_fractions}.
\\\\
\indent Both methods give very similar results, consistent with those obtained in \cite{saurel2008modelling, rodio2015innovative, pelanti2014mixture, zein2010modeling}. The density and volume fraction profiles obtained with the splitting approach show some overshooting in the middle region. This test case indicates that a strong rarefaction wave is well retrieved with both methods. Again, the newly proposed method takes larger time steps ($14303$ time steps) than the direct approach from \cite{daude2014numerical} ($14559$ time steps). \added[id=Rev.1]{The CPU time is $23.5$s and $27.4$s for the splitting approach and the direct approach, respectively (averaged over $10$ runs on an i5 processor).}

%
%

\section{Conclusions}\label{section Conclusions}
An acoustic-convective splitting\added[id=Authors]{-}based scheme has been proposed to solve the Kapila single-pressure single-velocity two-phase flow model. The acoustic and convective submodels are alternatingly stepped in time to approximate the solution of the entire flow model. The model dealing with the acoustic waves has been cast into a Lagrangian form, and solved using an HLLC-type solver. This approach gives a simple numerical scheme. The model dealing with the convective \added[id=Authors]{waves} has been approximated using a classical upwind scheme. The method has been evaluated for a variety of shock tube problems, and compared with an existing HLLC-type scheme applied to the original (unsplit) Kapila model. The obtained numerical results demonstrate the ability of the proposed method to deal with strong discontinuities and mixture flows. \added[id=Authors]{They} are in good agreement with exact and approximate reference solutions. The newly proposed method takes larger time steps than the HLLC-type scheme does for the unsplit model originally proposed in \cite{daude2014numerical}. \added[id=Rev.1]{This is most significant in the transonic regime.} Contact discontinuities, rarefaction waves and shock waves are captured very accurately with \added[id=Authors]{both the new method and the direct approach. The new method seems to be less diffusive than the direct approach.} Furthermore, the splitting approach may circumvent the inaccuracies when using approximate Godunov approaches for subsonic flows. The potential of the current method to deal with low-Mach number flows is briefly described in \cite{ten2015compressible}. \added[id=Revs]{To obtain higher-order temporal accuracy, the combination of higher-order methods to solve the systems together with a higher-order splitting approach must be used. One approach could be to use a generalized-$\alpha$ or a Runge-Kutta time integrator combined with a Strang splitting approach.} \added[id=Authors]{The proposed approach has a natural extension to multi-dimensional problems.}

%
%


\bibliographystyle{unsrt}
\bibliography{references}

\end{document}